\documentclass[12pt]{amsart}
\usepackage{fullpage}
\usepackage{amsmath,amsthm,amssymb,amscd}
\usepackage{graphicx,psfrag,epsfig} 

\begin{document}
\title{Generic diffeomorphisms with superexponential growth of number 
of periodic orbits.}
\author{Vadim Yu. Kaloshin}
\address{Department of Mathematics, 
Princeton University, Princeton NJ 08544-1000}
\email{kaloshin@math.princeton.edu}
\date{January 24, 1999}
\def\IMSmarkvadjust{0 pt}
\def\IMSmarkhadjust{0 pt}
\def\SBIMSMark#1#2#3{
 \font\SBF=cmss10 at 10 true pt
 \font\SBI=cmssi10 at 10 true pt
 \setbox0=\hbox{\SBF Stony Brook IMS Preprint \##1}
 \setbox2=\hbox to \wd0{\hfil \SBI #2}
 \setbox4=\hbox to \wd0{\hfil \SBI #3}
 \setbox6=\hbox to \wd0{\hss
             \vbox{\hsize=\wd0 \parskip=0pt \baselineskip=10 true pt
                   \copy0 \break%
                   \copy2 \break%
                   \copy4 \break}}
 \dimen0=\ht6   \advance\dimen0 by \vsize \advance\dimen0 by 8 true pt
                \advance\dimen0 by -\pagetotal
	        \advance\dimen0 by \IMSmarkvadjust
 \dimen2=\hsize \advance\dimen2 by .25 true in
	        \advance\dimen2 by \IMSmarkhadjust

%
%
  \openin2=publishd.tex
  \ifeof2\setbox0=\hbox to 0pt{}
  \else 
     \setbox0=\hbox to 3.1 true in{
                \vbox to \ht6{\hsize=3 true in \parskip=0pt  \noindent  
                {\SBI Published in modified form:}\hfil\break
                \input publishd.tex 
                \vfill}}
  \fi
  \closein2
  \ht0=0pt \dp0=0pt
 \ht6=0pt \dp6=0pt
 \setbox8=\vbox to \dimen0{\vfill \hbox to \dimen2{\copy0 \hss \copy6}}
 \ht8=0pt \dp8=0pt \wd8=0pt
 \copy8
 \message{*** Stony Brook IMS Preprint #1, #2. #3 ***}
}

\maketitle
\thispagestyle{empty}
\SBIMSMark{1999/2}{February 1999}{}

\markboth{Periodic orbits}{V. Kaloshin}       
 
\theoremstyle{plain}
\newtheorem{Thm}{Theorem}
\newtheorem{Def}{Definition}
\newtheorem{Lm}{Lemma}
\newtheorem{Prop}{Proposition}
\newtheorem{Corollary}{Corollary}
\newtheorem{Rem}{Remark}

\def\bdef{\begin{Def}}
\def\endef{\end{Def}}
\def\bthm{\begin{Thm}}
\def\ethm{\end{Thm}} 
\def\bprop{\begin{Prop}}
\def\enprop{\end{Prop}}
\def\bcor{\begin{Corollary}}
\def\encor{\end{Corollary}} 
\def\blm{\begin{Lm}}
\def\elm{\end{Lm}}
\def\brm{\begin{Rem}}
\def\erm{\end{Rem}}
\def\bfig{\begin{picture}}
\def\efig{\end{picture}}
\def\beq{\begin{eqnarray}}
\def\eneq{\end{eqnarray}}
\def\Cal{\mathcal}
\def\Bbb{\mathbb}
\def\R{\Bbb R}
\def\C{\Bbb C}
\def\Z{\Bbb Z}
\def\~{\tilde}
\def\dt{\delta}
\def\al{\alpha}
\def\eps{\epsilon}
\def\lb{\lambda}
\def\gm{\gamma}
\def\td{\tilde}

\hfill \begin{minipage}[t]{13em}
{Dedicated to the memory of my grandfather
Meyer Levich}
\end{minipage}

\section{Introduction}

Let $C^r(M,M)$ be the space of $C^r$ mappings of a compact 
manifold $M$ into itself with the uniform $C^r$-topology
and Diff$^r(M)$ be the space of $C^r$ 
diffeomorphisms of $M$ with the same topology. 
It is well-known that Diff$^r(M)$ is an open subset of
$C^r(M,M)$.
For a map $f \in C^r(M)$, consider the number of 
{\it{isolated}} periodic points of period $n$ (i.e. the
number of isolated fixed points of $f^n$)
\beq \label{grow}
P_n(f)=\# \{\ {\textup{isolated}}\ \ x \in M:  x=f^n(x)\}.
\eneq 
In 1965 Artin \& Mazur  {\cite {AM}} showed
that: {\it there exists a dense set $\Cal D$ in $C^r(M,M)$
such that for any map $f \in \Cal D$ the number $P_n(f)$
grows at most exponentially with $n$, i.e. for some number $C>0$} 
\beq
P_n(f) \leq \exp(Cn) \ \  {\textup{for all}}\ \  n \in \Bbb Z_+.
\eneq

Notice that the Artin-Mazur Theorem does not exclude the 
possibility that a mapping $f$ in $\Cal D$ has a curve of periodic 
points $\gm$, i.e. $\forall x \in \gm,\ f^n(x)=x$ for some 
$n\in\Z_+$, because in this case $\gm$ consists of nonisolated
periodic points of period $n$ (see the last part of Theorem \ref{mthm} 
for this nonisolated case).

\bdef
We call a mapping (resp. diffeomorphism)  $f \in C^r(M,M)$
(resp. $f \in \textup{Diff}^r(M)$)
an {\it Artin-Mazur 
mapping (resp. diffeomorphism) or simply A-M 
mapping (resp. diffeomorphism)} if
$P_n(f)$ grows  at most exponentially fast. 
\endef

Artin-Mazur {\cite{AM}} posed  the following 
problem: {\it{what can be said about the set of A-M mappings with only
transversal periodic orbits  in the space of $C^r$ mappings?}}
Recall that a periodic orbit of period $n$ is called 
{\it{transversal}} if the linearization $df^n$ at this point 
has for an eigenvalue no $n$-th roots of unity. Notice that a 
hyperbolic periodic point is always transversal, but not vice versa. 

In what follows we consider {\it not the whole space $C^r(M,M)$ 
of mappings of $M$ into itself, but only its open subset 
{\textup{Diff}}$^r(M)$.} 
The first main result of this paper is an 
answer to Artin-Mazur's question for the space of 
diffeomorphisms and a simple proof of it.

\bthm \label{AMthm}
Let $1 \leq r < \infty$. Then the set of A-M diffeomorphisms 
with only hyperbolic periodic orbits is dense in the space 
{\textup{Diff}}$^r(M)$ with the uniform $C^r$ topology. 
\ethm

In a weak form this theorem says that A-M diffeomorphisms which
satisfy the Kupka-Smale condition form a dense set in
${\textup{Diff}}^r(M)$. Recall that a diffeomorphism 
is called {\it{a Kupka-Smale (or K-S) diffeomorphism)}} if all its 
periodic points are hyperbolic and all associated stable and unstable
manifolds intersect one another transversally. The Kupka-Smale theorem 
says that K-S diffeomorphisms form a residual set (see e.g. {\cite{PM}}). 
The natural question is whether intersection of A-M and K-S 
diffeomorphisms can be dense in Diff$^r(M)$. The answer is not easy,
because methods of the proof of both theorems are of completely 
different nature and can not be applied simultaneously.
If one omits the condition on transversality of stable and unstable manifolds, 
then Theorem \ref{AMthm} says that the intersection of A-M
and K-S diffeomorphisms is dense. 

A residual set in a finite-dimensional space can have Lebesgue
measure zero. Therefore, the Kupka-Smale theorem does not imply 
that ``almost every'' diffeomorphism is a K-S diffeomorphism.
In loose terms, a set $P \subset {\textup{Diff}}^r(M)$ is 
called {\it{prevalent}} (or ``has a full measure'') if for a generic 
finite-parameter family $\{f_\eps\}_{\eps \in \textup{Ball}}$, 
the property $f_\eps \in P$ holds for almost every parameter value. 
In {\cite{K}} it is proven that K-S diffeomorphisms form a prevalent set
or have ``a full measure''. Finally, for a discussion of prevalence see 
section \ref{prev}.

In {\cite{AM}} Artin-Mazur also introduced the {\it {dynamical 
$\zeta_f$-function}} defined by
$$
\zeta_f(z)=\exp\left(\sum_{n=1}^\infty P_n(f) \frac{z^n}{n}\right).
$$

For A-M diffeomorphisms the dynamical $\zeta_f$-function is
analytic in some disk centered at zero. It is well-known that the
dynamical $\zeta_f$-function of a diffeomorphism $f$ satisfying Axiom A 
has an analytic continuation to a rational function (e.g. {\cite{Ba}}). 

In 1967 Smale {\cite {S}} posed the following
question (Problem 4.5, p.765):

{\it {Is the dynamical $\zeta_f$-function generically rational
(i.e. is $\zeta_f$ rational for a residual set of 
$f \in {\textup{Diff}}^r(M)$)?}}

In {\cite {Si} it is shown that for the $3$-dimensional 
torus the  $\zeta_f$-function is not rational.
It turns out  that for manifolds of dimension greater
or equal than $2$ it is not even analytic in any
neighborhood of zero (see Theorem \ref{mthm} below).
Recall that a subset of a topological space is called residual 
if it contains a countable intersection of open dense subsets.
We call a residual set a {\it{topologically generic}} set.

Finally, in 1978 R. Bowen asked the following question in his book {\cite{Bo}}:

Let $h(f)$ denote the  topological entropy of $f$.
{\it{Does }} 
$$h(f)=\limsup_{n \to \infty} {\log P_n(f)}/{n}$$
{\it{for a generic set of diffeomorphisms with respect to the
$C^r$ topology?}}

It turns out the two above questions can be answered simultaneously
for $C^r$ diffeomorphisms with $2 \leq r < \infty$. 
The second result is the following:
\bthm \label{main} Let $2 \leq r < \infty$. Then
the set of A-M diffeomorphisms is {\bf not} topologically generic 
in the space of $C^r$ diffeomorphisms
{\textup{Diff}}$^r(M)$ with the uniform $C^r$ topology.
\ethm
We have the following consequences:
\bcor The property of having a convergent $\zeta_f(z)$
function is not $C^r$-generic, nor is the equation
$h(f)=\limsup_{n \to \infty} {\log P_n(f)}/{n}$.
\encor

The first part is easy. To prove the second,
notice that the topological entropy for any $C^r$ ($r \geq 1$)
diffeomorphism $f$ of a compact manifold is always finite 
(see e.g. {\cite{HK}}). Define
the rate of growth of the number of periodic orbits  by
$\limsup_{n \to \infty} {\log P_n(f)}/{n}$.
Then for diffeomorphisms which are not A-M diffeomorphisms,
the rate of growth is always equal to infinity.

Since, an Axiom A diffeomorphism is an A-M diffeomorphism, we
need to analyze the complement to Axiom A diffeomorphisms in the space
of {\textup{Diff}}$^r(M)$.
An example of a diffeomorphism with an arbitrarily fast 
growing number of periodic orbits is given in {\cite {RG}}.  
Now we describe a ``bad'' domain, where the A-M property fails to be
topologically generic. 

In 1970 Newhouse found a domain in the space of $C^r$ 
diffeomorphisms Diff$^r(M)$, where
diffeomorphisms exhibiting homoclinic tangencies are dense \cite{N}.
Such a domain is called {\it a Newhouse domain} 
$\Cal N \subset {\textup{Diff}}^r(M)$. Our third and main result is as follows:

\bthm \label{mthm} Let $2 \leq r <\infty$. Let   
$\Cal N \subset {\textup{Diff}}^r(M)$ be a Newhouse domain. Then 
for an arbitrary sequence of positive integers        
$\{ a_n \}_{n=1}^\infty$ there exists a 
residual set $\Cal R_a \subset \Cal N$, depending on 
the sequence $\{ a_n \}_{n=1}^\infty$,
with the property that $f \in \Cal R_a$ implies that
$$
\limsup_{n \to \infty} 
{P_n(f)}/{a_n}=\infty.
$$ 
Moreover, there is a dense set $\Cal D$ in $\Cal N$ such that
any diffeomorphism $f \in F$ has a curve of periodic points.
\ethm

Let us note that Theorems  {\ref{main}} and {\ref{mthm}} 
follow from a Theorem of 
Gonchenko-Shilnikov-Turaev which will be discussed in section 2.

In such a domain Newhouse exhibited a residual set of 
diffeomorphisms with infinitely many distinct sinks 
\cite{N}, \cite{R}, and \cite{PT}.  
Now it is known as {\it Newhouse's phenomenon}.
In a way Theorem \ref{mthm} is similar to Newhouse's phenomenon
in the sense that for a residual set a ``bad'' property holds true.

Continuing, Theorem \ref{main} is a corollary of the first part of 
Theorem \ref{mthm}. To see this 
fix the sequence $a_n=n^n$ and denote by $\Cal R_a$ a set from Theorem 
\ref{mthm} corresponding to this sequence. Assume that A-M 
diffeomorphisms form a residual set, then this set must intersect 
with $\Cal R_a$ which is a contradiction.

It seems that based on Newhouse's phenomenon in the space 
Diff$^1(M)$ with the $C^1$-topology, where $\dim M \geq 3$, found by 
Bonnati \& Diaz {\cite{BD}} one can extend Theorems {\ref{main}} 
and {\ref{mthm}} to the case $r=1$ and $\dim M \geq 3$.
The problem with this straightforward generalization is
that the proof of the Gonchenko-Shilnikov-Turaev  (GST) theorem
is essentially two-dimensional. To generalize the GST theorem to 
the three-dimensional case one needs either to find an invariant
two-dimensional surface and use the two-dimensional proof or 
find another proof. In personal communications, Lorenzo 
Diaz has shown to the author that an invariant two-dimensional 
surface can be constructed using the method from {\cite{BD}}.
However, this extension is not  straightforward and will appear
separately from this paper. 

Analogs of Theorems  {\ref{main}}  and {\ref{mthm}} can be formulated
for the case of vector fields on a compact manifold of dimension 
at least $3$. Reduction from the 
case of diffeomorphisms to the case of vector fields can be
done using the standard suspension of a vector field over 
a diffeomorphism {\cite {PM}}.

\subsection{Newhouse phenomenon and Palis conjecture}
Newhouse showed that a Newhouse domain exists under
the following hypothesis:

Let a diffeomorphism $f \in {\textup{Diff}}^r(M)$ have a saddle 
periodic orbit $p$. Suppose stable $W^s(p)$ and unstable $W^u(p)$
manifolds of $p$ have a quadratic tangency. Such a
diffeomorphism $f$ is called                       
{\it{a diffeomorphism exhibiting a homoclinic tangency}}.
Then arbitrarily 
$C^r$-close to $f$ in Diff$^r(M)$ there exists {\it a Newhouse domain}. 
In particular, it means that by a small 
$C^r$-perturbation of a diffeomorphism $f$ with a homoclinic 
tangency one can generate arbitrarily quick growth of the
number of periodic orbits. 

On this account we would also like to mention the following
conjecture, which is due to Palis {\cite {PT}}, about the space
of diffeomorphisms of $2$-dimensional manifolds:

{\bf {Conjecture}} {\it {If $\dim M=2$, then every 
diffeomorphism $f \in {\textup{Diff}}^r(M)$ can be approximated by
a diffeomorphism which is either hyperbolic or 
exhibits a homoclinic tangency}}.

This conjecture is proven for approximations in the $C^1$ topology 
{\cite{PS}}. If this conjecture is true, 
then in the complement to the set of hyperbolic 
diffeomorphisms those diffeomorphisms with arbitrarily quick growth of 
number of periodic orbits form a topologically generic set.

Unfolding of homoclinic tangencies is far from being understood.
In {\cite{GST}} the authors describe the following important result:
{\it{there does not exist a finite number of parameters to describe
all bifurcations occurring next to a
homoclinic tangency}} (see section 2, Corollary \ref{cthree} for details). 
This implies that the complete description of
bifurcations of diffeomorphisms with a homoclinic tangency 
is impossible.

This paper is organized as follows. In section 2.1 we state the 
Gonchenko-Shilnikov-Turaev result and give a proof of it
in sections 2.2--2.5. 
Section 2.6 is devoted to the proof
of Theorem {\ref{mthm}} in the case $\dim M=2$. 
Reduction from Theorem {\ref{mthm}} in the case $\dim M=2$ to 
Theorem {\ref{mthm}} to the general case $\dim M\geq 2$ is done in
section 2.7. Theorem \ref{AMthm} is proven in part 3.
Discussion of open questions related to the main results is in section 4. 

From now on we consider
diffeomorphisms of a 2-dimensional compact manifold $M$.
The general case when $\dim M \geq 2$ can be reduced to the 
$2$-dimensional case by the standard suspension and the Fenichel theorem 
{\cite {Fe}} (see section \ref{multi}). 

\section{Degenerate periodic orbits in a Newhouse domain and 
the Gonchenko-Shilnikov-Turaev Theorem {\cite {GST}}}

Assume that a $C^r$ diffeomorphism $f$ exhibits a
homoclinic tangency. By the Newhouse theorem {\cite {N}}, in each 
$C^r$ neighborhood  of a diffeomorphism
$f$ exhibiting a homoclinic tangency there exists a Newhouse domain.

Let us define a degenerate periodic point of order $k$ or
a $k$-degenerate periodic point. Sometimes, it is also
called a {\it{saddlenode periodic orbit of multiplicity}} $k+1$. 
 
\bdef Let $f$ be a $C^s$ diffeomorphism of a $2$-dimensional 
manifold having a periodic orbit $p$ of period $m$. A
periodic point $p$ is called $k$-degenerate, where $k < s$, if 
the linear part of $f^m$ at point $p$ has a multiplier 
$\nu =1$ while the other multiplier is different in 
absolute value from the unit and a restriction of $f$ to 
the central manifold in some coordinate system can be written 
in the form
\beq \label{deg}
x \mapsto x+l_{k+1} x^{k+1} + o(x^{k+1}).
\eneq
\endef

Let $s>r$. Then $C^s$ diffeomorphisms are dense
in the space Diff$^r(M)$ and, therefore, in any 
Newhouse domain $\Cal N \subset {\textup{Diff}}^r(M)$ 
(see e.g. {\cite{PM}}).

\bthm (Theorem 4, {\cite {GST}}) {\label {Shil}}
For any positive integers $s > k \geq r$ 
the set of $C^s$ diffeomorphisms having a $k$-degenerate
periodic orbit is dense in a
Newhouse domain $\Cal N \subset {\textup{Diff}}^r(M)$.
\ethm

This theorem and Newhouse's theorem imply the following important result:

\bcor {\cite {GST}}\label{cthree} Let $f \in {\textup{Diff}}^r(M)$ be
a diffeomorphism exhibiting a homoclinic tangency. 
There is no finite number $s$ such that 
a generic $s$-parameter family $\{f_\varepsilon\}$ unfolding a
diffeomorphism $f_0=f$ is a versal 
family of $f_0$ meaning that the family $\{f_\varepsilon\}$
describes all possible bifurcations occurring 
next to $f$. Indeed, to describe all possible bifurcations of 
a $k$-degenerate periodic orbit one needs at least $k+1$ parameters and
$k$ can be arbitrary large.
\encor

Once Theorem \ref{Shil} is proved the proof of Theorem 
{\ref{mthm}} can be completed by inductive application of the 
following idea. Let $f$ be a $C^s$ diffeomorphism  from a Newhouse 
domain $\Cal N \subset {\textup{Diff}}^r(M)$ with a $k$-degenerate 
periodic orbit $p$ of period, say $n$, of $f$ for $s>k\geq r$, then $p$ is 
{\it{flat periodic point}} along the central manifold with respect 
to the $C^r$ topology, namely, by a $C^r$-perturbation one can 
make the restriction to the central manifold be {\it{the identical map}}.
It allows us either to create a {\it{curve of periodic orbits}} or
{\it{split}} $p$ into {\it{any ahead
given number of}} hyperbolic periodic orbits of the same period (or
double the period of $p$) by a small perturbation. Since,
created periodic orbits are hyperbolic they persist under
 perturbations. Moreover, after a perturbation we are
still in a Newhouse domain one can iterate this procedure 
of creating a $k$-degenerate periodic orbits and splitting them without 
destroying what was done in previous stages (see section \ref{growth}).

In what follows we need a few notions
related to a saddle periodic point. These definitions will be 
needed in the proof of Theorem {\ref {Shil}}.

\bdef Let $f$ be a $C^s$ diffeomorphism of\  a  
$2$-dimensional manifold $M$ and let $p$ be a saddle periodic point 
of period $m$, namely, $f^m(p)=p$ with eigenvalues $\lb$ and 
$\mu$, $\lb<1<\mu$. The saddle exponent of $p$ is the number
$\rho(p,f)=\frac{-\log \lb}{\log \mu}$. We call $p$  
a $\rho$-shrinking saddle, where $\rho=\rho(p,f)$. If $\rho$ is greater 
than some $r$, then $p$ is also called at least $r$-shrinking.

A saddle $p$ is called nonresonant if for any pair of 
positive integers $n$ and $m$ such that the
number $\lb^n \mu^m$ is different from $1$.
\endef

\subsection{A Scheme of a Proof of Theorem \ref{Shil}}

Theorem \ref{Shil} is stated in ({\cite{GST}}, Thm.4). 
A proof of this theorem is outlined there. 
Proof of several technical statements
\footnote{Lemmas 1 and 2 in {\cite{GST}} which corresponds to Lemmas \ref{normf} and \ref{lthree} of the present paper respectively } 
is omitted there. 
We present a rigorous proof which essentially uses ideas given in 
{\cite{GST}}. In what follows a $C^r$-perturbation means a small 
$C^r$-perturbation. The proof of Theorem \ref{Shil} consists of four steps.

{\it{The first step}}. 
From the existence of a homoclinic tangency of a dissipative saddle, we
deduce the existence (after a $C^r$-perturbation) of a 
homoclinic tangency of an at least $k$-shrinking saddle, $k>r$.

{\it{The second step}}. From the existence of a homoclinic tangency 
of an at least  $k$-shrinking saddle, we create a $k$-floor tower 
(defined in section 2.4) after a $C^r$-perturbation (see Fig.3 for $k=3$).

{\it{The third step}}. From the existence of
a $k$-floor tower, we show  that a $C^r$-perturbation 
can make a $k$-th order homoclinic tangency.

{\it{The fourth step}}. From the existence of a $k$-th order 
homoclinic tangency we construct by a $C^r$-perturbation
a $k$-th order degenerate periodic orbit of an arbitrarily high
period.

Notice that the way we construct a $k$-tower is slightly different 
from the one in {\cite{GST}}.

The proof of Theorem \ref{Shil}  is given in sections 2.2--2.5
according to the following plan.
In section 2.2 we present some basic properties of a return map in a 
neighborhood of a quadratic homoclinic tangency. In section 2.3
we realize the first step (Corollary \ref{step1}) and calculate 
limits for return maps in a neighborhood of a $k$-th order homoclinic 
tangency, where $k \geq 2$. The second and the third steps are done 
in sections 2.4 and 2.5 respectively. The last fourth step
consists in application of Corollary \ref{ktang} proven in section 2.3. 

\subsection{Basic properties of a return map in a neighborhood of
a homoclinic tangency}\label{s2.2} 

Fix a positive integer $r \geq 2$.
Consider a $C^\infty$ smooth diffeomorphism $f: M^2 \to M^2$ with a  
saddle fixed point $p$, namely, $f(p)=p$ with the eigenvalues $\lb$ and $\mu$. 
Assume the saddle $p$ is dissipative and nonresonant. 
We can obtain all conditions  by applying a $C^r$-perturbation 
(for $f \in C^r$, or/and $\lb \mu=1$, or/and by inverting $f$) if necessary. 
Then by the standard fact from the theory of normal forms e.g. {\cite{IY}} 
the map $f$ is $C^{r}$ linearizable in a neighborhood $U$ of $p$
\beq \label{lin}
f:(x,y) \mapsto (\lambda x,\ \mu y),
\eneq 
where $\lambda < 1 < \mu$ and $\lambda \mu <1$. The larger is $r$, the
smaller is the neighborhood $U$, where a $C^r$-normal form applicable.

Assume that the stable $W^s(p)$ and unstable manifold $W^u(p)$ of 
$p$ in normal coordinates have a point of quadratic tangency 
$q$ with coordinates $(1,0)$ and for some $N$ we have that 
$f^{-N}(q)=\tilde q$ has coordinates $(0,1)$ (see Fig. 1).
Assume also that  in a neighborhood of the  
homoclinic point $q$ the unstable manifold $W^u(p)$ lies in the 
upper half plane $\{y \geq 0\}$ and the directions of $W^u(p)$ and
$W^s(p)$ at the point of tangency $q$ are the same (see Fig.1).
Diffeomorphisms with such type of homoclinic tangency are dense 
in a Newhouse domain see e.g. \cite{PT}.
 
\begin{figure}[htbp]
  \begin{center}
   \begin{psfrags}
     \psfrag{(0,1)}{\small{$(0,1)$}}
     \psfrag{(1,0)}{\small{$(1,0)$}}
     \psfrag{WUP}{\small{$W^u(p)$}}
     \psfrag{WSP}{\small{$W^s(p)$}}
     \psfrag{tilq}{\small{$\tilde q$}}
     \psfrag{q}{\small{$q$}}
     \psfrag{p}{\small{$p$}}
    \includegraphics[width= 2in,angle=0]{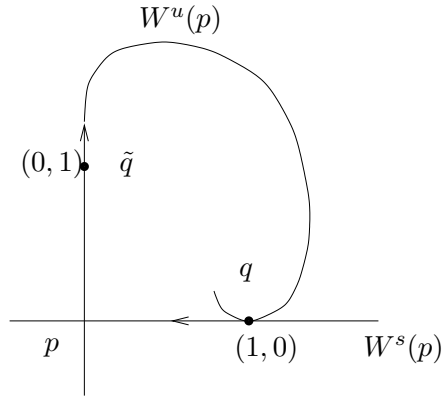}
    \end{psfrags}
   \caption{Homoclinic tangency}
  \end{center}
\end{figure}

Rescale sufficiently small neighborhoods 
$U$ and $\tilde U$ of $q$ and $\tilde q$
respectively. Denote by $W_{loc}^u(p)$ the first connected component
of the intersection $W^u(p) \cap U$ (see Fig. 2). Below we shall
use the coordinate systems in $U$ and $\tilde U$ induced by the 
normal coordinates of $p$ and $f$. Write $W^u_{loc}(p)$ in $U$ as the 
graph of a function $y=cx^2+g(x)$, where $g(x)=o(x^2),\ c>0$.
A rectangle in $U$ (resp. $\tilde U$)
is called {\it{a right rectangle}} if it has two side that are parallel 
to the coordinate axis. 
 
\begin{figure}[htbp]
  \begin{center}
    \begin{psfrags}
      \psfrag{(0,1)}{\footnotesize{$(0,1)$}}
      \psfrag{(1,0)}{\footnotesize{$(1,0)$}}
      \psfrag{ln}{\small{$l_n$}}
      \psfrag{dtn}{\small{$\dt_n$}}
      \psfrag{Dn}{\small{$\Delta_n$}}
      \psfrag{dtnr}{$\dt_n^\rho$}
      \psfrag{Ln}{\small{$L_n$}}
      \psfrag{tn}{\small{$\tau_n$}}     
      \psfrag{WUP}{\small{$W^u(p)$}}
      \psfrag{WSP}{\small{$W^s(p)$}}
      \psfrag{app}{\small{$\sim$}}
     \includegraphics[width= 2.8in,angle=0]{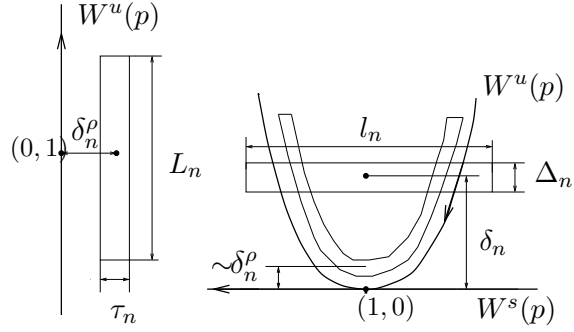}
     \end{psfrags}
   \caption{Neighborhoods of points of homoclinic tangencies.}      
  \end{center}
\end{figure}

\bprop \label{return} 
Let $c>0$ be the above constant  and $n \in \Bbb Z_+$
be sufficiently. 
Put $\delta_n=\mu^{-n}, \Delta_n=2\delta_n^{3/2}$,
and $l_n=3c\delta_n^{1/2}$. Consider a right 
rectangle $T_n$ centered at $(1,\dt_n)$ whose horizontal
(resp. vertical) length is $l_n$ (resp. $\Delta_n$). 
Then the image $f^n(T_n)$ is the right rectangle centered at 
$(\dt_n^\rho,1)$ whose horizontal
(resp. vertical) length is $\tau_n=3c \delta_n^{1/2} \lambda^n$ 
(resp. $L_n= 2\delta_n^{1/2}$).
Moreover, $T_n$ and $f^{n+N}(T_n)$ form a horseshoe 
which has a periodic saddle $p'$ of period $n+N$ and 
the curvilinear rectangle $f^{n+N}(T_n)$ is
$C\delta_n^\rho$ distant away from $W^s(p)$ for some $C>0$ independent 
of $n$ (see Fig. 2). 
\enprop

\brm
The exponent $\rho=\frac{\log 1/\lb}{\log \mu}$ gives a 
characteristic of distortion while a trajectory pass in 
a neighborhood of saddle $p$ in the following sense. 
The rectangle $T_n$ is $\dt_n$-distance away from $W^s(p)$, but 
its image $f^{n+N}(T_n)$ is $\dt_n^\rho$-distance away from $W^s(p)$, so
the more exponent $\rho=\frac{\log 1/\lb}{\log \mu}$
of saddle $p$, the deeper the horseshoe part  $f^{n+N}(T_n)$ 
goes inside $W^u_{loc}(p)$ with respect to $T_n$ and gets
closer to $W^s(p)$. 
\erm
{\it Proof of Proposition \ref{return}:}\ \ 
Use the fact that $f$ is linear (\ref{lin}) in
the unit square around $p$ to prove the first part.
Since $N$ is fixed the ratio of \{distance of the rectangle $f^n(T_n)$
to $W^u(x)$\} and \{distance of the curvilinear rectangle $f^{n+N}(T_n)$
to $W^u_{loc}(x)$\} is bounded. This implies the second statement
of the Proposition and completes the proof. Q.E.D.

\subsection{The first step: higher order homoclinic tangencies and appropriate limits}

It is well-known from e.g. {\cite {MV}}, {\cite {PT}}, and {\cite {TY}}
that for a generic 1-parameter family $\{f_\eps\}$ unfolding 
a quadratic (1-st order) homoclinic tangency $q$ and for any number 
$c \in [-2,1/4]$ there exist three sequences
($n \in \Bbb Z_+$): rectangles $T_n$ next to $q$, rescalings 
$R_n:T_n \to [-2,2]\times [-2,2]$, and parameters $\eps_n$ 
such that a return map $f_{\eps_n}^n$ from $T_n$ into itself converges to the map 
$(x,y) \to (y,y^2+c)$. $T_n$ and $R_n$ are  
independent of $c$, but $\eps_n=\eps_n(c)$ depends on $c$.

In this section we show that for a generic $k$-parameter family
$\{f_\mu\}$ unfolding a $(k-1)$-st order homoclinic tangency $q$ and 
for any set of numbers $M=(M_0, \dots ,M_{k-1}) \in \Bbb R^k$
there exist three sequences ($n \in \Bbb Z_+$):
rectangles $T_n$ next to $q$, rescalings 
$R_n:T_n \to [-2,2]\times [-2,2]$, and parameters 
$\mu(n)=(\mu_0(n),\dots, \mu_{k-1}(n))$ such  that a return map $f_{\mu(n)}^n$
from $T_n$ into itself under $f_{\eps_n}$ converges to 
$(x,y) \to (y,y^k+\sum_{i=0}^{k-1}M_i y^i).$
Moreover, convergence holds with respect to the $C^r$ topology 
for any $r \in \Bbb Z_+$. The calculation presented below is 
in the spirit of {\cite {PT}} and {\cite {TY}}.

Consider a $C^\infty$ diffeomorphism $f$ which has a dissipative saddle
periodic point $p$ exhibiting a homoclinic tangency of $(k-1)$-st order. 
Figure 1  illustrates the topology for even $k$.
We shall use notations of Fig. 1.
Take coordinates $(\bar x,\bar y)=(x-1,y)$
around the homoclinic point $q$ and coordinates 
$(\tilde x,\tilde y)=(x,y-1)$
around the homoclinic point $\tilde q$. 

Because of a $(k-1)$-st order tangency, the map 
$f^N:(\tilde x,\tilde y) \to (\bar x,\bar y)$ from 
a neighborhood $\tilde U$ of $\tilde q$  with coordinates 
$(\tilde x,\tilde y)$ to a neighborhood $U$ of 
$q$ with coordinates $(\bar x,\bar y)$ can be written in the form:

\beq \label{outmap}
\left( \begin{array}{c}
\tilde x \\
\tilde y
\end{array} \right)
\stackrel{f^N}{\longmapsto} 
\left( \begin{array}{c}
\alpha \tilde y+\ \beta \tilde x +\ H_1(\tilde x,\tilde y) \\
\gamma \tilde y^{k} +\ \sigma \tilde x+\ H_2(\tilde x,\tilde y)
\end{array} \right),
\eneq
where $\alpha,\ \beta,$ and $\gamma$ are constants such that
for $\tilde x=\tilde y=0$

\beq \label{remainder}
\begin{cases}
H_1=\partial_x H_1=\partial_y H_1=0 \\ 
H_2=\partial_x H_2=\partial_y^j H_2=0,\ \ j=1,\dots, k.
\end{cases} 
\eneq
To see that formula (\ref{outmap}) holds consider images 
of lines $\{\tilde x={\textup{constant}}\}$.

Consider a generic $k$-parameter unfolding of a $(k-1)$-st order 
homoclinic tangency:

\beq \label{unfold}
\left( \begin{array}{c}
\tilde x \\
\tilde y
\end{array} \right)
\stackrel{f^N_\mu}{\longmapsto} 
\left( \begin{array}{c}
\alpha \tilde y+\ \beta \tilde x +\ H_1(\tilde x,\tilde y) \\
\gamma \tilde y^{k} +\ \sum_{i=0}^{k-1} \mu_i \tilde y^i+\ 
\sigma \tilde x+\ H_2(\tilde x,\tilde y)
\end{array} \right).
\eneq

The main result of this section in the following:

\blm \label{normf} With the above notations and $k \geq 2$
for an arbitrary set of real numbers $\{M_i\}_{i=0}^{k-1}$
there exists a sequence of parameters $\{\mu(n)\}_{n \in \Bbb Z_+}$
such that $\mu(n)$ tends $0$ as $n \to \infty$ 
\footnote{see (\ref{par}) for the exact form of the sequence} 
and a sequence of change of variables 
$R_n:(\bar x, \bar y) \to (x_n,y_n)$ such that 
the sequence of maps:
$\{R_n \circ f^{n+N}_{\mu(n)} \circ R_n^{-1}:[-2,2]\times [-2,2] \to 
[-2,2]\times [-2,2]\}$ converges to the 1-dimensional map 
\beq \label{normform}
\left( \begin{array}{c}
x \\ y
\end{array} \right)
\stackrel{\phi_M}{\longmapsto} 
\left( \begin{array}{c}
y \\ y^{k}+\sum_{i=0}^{k-1} M_i y^i 
\end{array} \right)
\eneq 
in the $C^r$ topology for any $r$.
\elm
\bcor \label{ktang} (The fourth step of the proof of Theorem \ref{Shil}) 
For $M_1=1$, and $M_j=0,\ j=0,2, \dots, k-1$    
by a $C^r$-perturbation of a $C^\infty$ diffeomorphism $f$ exhibiting a 
homoclinic tangency of order $(k-1)$
one can create a $C^r$ diffeomorphism $\tilde f$
with a $(k-1)$-degenerate periodic orbit of an arbitrarily high
period.  
\encor

\bcor \label{step1}
For $k=2$, $M_0=-2$, and $M_1=0$ by a $C^r$-perturbation 
of a $C^\infty$ diffeomorphism 
$f$ exhibiting  a quadratic homoclinic tangency
one can create a $C^\infty$ diffeomorphism 
$f$ with a periodic saddle $p$ exhibiting a homoclinic tangency and 
eigenvalues of $p$ are close to $2$ and to $+0$ respectively. 
\encor
To prove this corollary recall that for any $r$ the map
$(x,y) \to (y,y^2-2),$ $x,y \in [-2,2]$ has a fixed point $(2,2)$.
One can show that 
by a $C^r$ perturbation of this 2-dimensional map a fixed point $(2,2)$ 
becomes a saddle near to  $(2,2)$ exhibiting a homoclinic tangency.
In {\cite {PT}} \S6.3 Prop.3 figures 6.4 and 6.5 or 
{\cite {MV}} pg.14, this is shown to be true.
On 2-dimensional perturbations of the 1-dimensional map 
$y \mapsto y^2-\mu$ see also {\cite {BC}}.

{\it{Proof of Lemma}} \ref{normf}:
We follow the standard method  
and split the return map $f^{n+N}$ into the composition of
two maps: the linear map $f^n:(x,y) \to (\lambda^n x, \mu^n y)$ 
and  the map $f^N_\mu$ given by formula (\ref{unfold}).
The composition of $f^N_\mu$ and $f^n$  has the form:

\beq \label{compos}
\left( \begin{array}{c}
\bar x \\
\bar y
\end{array} \right)
\stackrel{f^N_\mu \circ f^n}{\longmapsto}
\left( \begin{array}{c}
\alpha \bar y_n+\ \beta \lambda^n(1+\bar x) +\ 
H_1(\cdot,\cdot)\\
\gamma \bar y_n^{k} +\sum_{i=0}^{k-1} \mu_i \bar y_n^i+
\sigma \lambda^n(1+\bar x)+H_2(\cdot,\cdot) 
\end{array} \right).
\eneq
where $\bar y_n=\ \mu^n\bar y-1$, 
$H_j(\cdot,\cdot)=\ H_j(\lambda^n(1+\bar x),\bar y_n)$, $j=1,2$.
Denote $\mu^{1/(k-1)}$ by $\tau$. Introduce the change of variables
$R_n:(\bar x,\bar y) \to (x_n,y_n)$, where
\beq
\left( \begin{array}{c}
x_n \\
y_n
\end{array} \right)
=
\left( \begin{array}{c}
\tau^n \bar x \\ 
\tau^{n} (\mu^n\bar y-1) 
\end{array} \right).
\eneq

In $(x_n,y_n)$-coordinates the map $f^N_\mu \circ f^n$
has the form:

\beq \nonumber
\left( \begin{array}{c}
x_n \\
y_n
\end{array} \right)
\stackrel{f^N_\mu \circ f^n}{\longmapsto} 
\left( \begin{array}{c}
\alpha y_n+\ \beta \lambda^n (\tau^n +\ x_n)
+\ \tau^{n} H_1(\cdot,\cdot) \\
\gamma y_n^{k} +\sum_{i=0}^{k-1} \mu_i \tau^{n(k-i)}y_n^i+
\sigma \lambda^n \mu^n(\tau^n+x_n)+
\tau^{kn} H_2(\cdot,\cdot) -\tau^n
\end{array} \right),
\eneq
where $H_j(\cdot,\cdot)=H_j(\lambda^n (1+\tau^{-n} x_n),\tau^{-n} y_n)$ 
for $j=1,2$.

Recall that $p$ is dissipative, so $\lb \mu<1$ and $\lb \tau<1$ too.
Thus, condition (\ref{remainder}) and $0<\lambda,\  \tau^{-1} <1$
imply that terms $\lambda^n \mu^n x_n,\ \beta \lambda^n (\tau^n +\ x_n),\ 
\tau^{n}H_1(\lambda^n (1+\tau^{-n} x_n),\tau^{-n} y_n)$,
and $\tau^{kn} H_2(\lambda^n (1+ \tau^{-n} x_n),\tau^{-n} y_n)$
tends to $0$ as $n \to \infty$ in the $C^r$ topology 
for any positive integer $r$.
 
Put
\beq 
\begin{aligned}\label{par} 
\mu_0(n)=\mu^{-kn/(k-1)}M_0-\sigma \lb^{n}+\mu^{-n} \\
\mu_i(n)=\mu^{-(k-i)n/(k-1)}M_i \ \ {\textup{for}}\ \  i=1, \dots, k-1.
\end{aligned}
\eneq
We see that all $\{\mu_i(n)\}$ tends to $0$ as $n$ tends to infinity.
Therefore, in the limit as $n \to \infty$ we obtain 
\beq \label{unfolding}
\left( \begin{array}{c}
\tilde x \\
\tilde y
\end{array} \right)
\stackrel{\phi_M }{\longmapsto} 
\left( \begin{array}{c}
\alpha y \\
\gamma y_n^{k} +\ \sum_{i=0}^{k-1} M_i y^i
\end{array} \right).
\eneq
Additional change of variables depending on $\alpha$ and 
$\gamma$ completes the proof. Q.E.D.

\subsection{The second step: Construction of a $k$-floor tower}

Consider a $C^\infty$ diffeomorphism $f$ with a nonresonant saddle
periodic point $p$ exhibiting a homoclinic tangency at a point $q$. 
First, we give a definition of {\it{a $k$-floor tower}}.
Recall that $U$ denotes a neighborhood of the homoclinic tangency
$q$. Let $\tilde p$ be a saddle periodic orbit of $f$,\ $\tilde p \in U$.
Then denote by $W^s_{loc}(\tilde p)$ (resp. $W^u_{loc}(\tilde p)$) the first 
connected component of the intersection of stable (resp. unstable)
manifold $W^s(\tilde p)$ (resp. $W^u(\tilde p)$) with $U$.

\bdef A $k$-floor tower is a contour consisting of $k$ saddle
periodic points $p_1, \dots ,p_r$ (of different periods) 
such that $W^u_{loc}(p_i)$ is tangent to $W^s_{loc}(p_{i+1})$
for $i=1,\dots , k-1$, and $W^u_{loc}(p_k)$ intersects
$W^s_{loc}(p_{1})$ transversally (see Fig.3 for $k=3$).
\endef

Construction of a $k$-floor tower is an intermediate step
in the proof of Theorem {\ref{Shil}}.
In this section we prove that
 
\blm For any positive integer $k$ a $C^r$ diffeomorphism $f$ 
exhibiting a homoclinic tangency for an at least $r$-shrinking saddle 
periodic orbit $p$ admits a $C^r$-perturbation $\tilde f$ such that
 $\tilde f$ has a $k$-floor tower. If $q$ is
a point of homoclinic tangency of $f$, then the aforementioned
tower of $\tilde f$ is located in a neighborhood $U$ of 
$q$.
\elm
 
\begin{figure}[htbp]
  \begin{center}
   \begin{psfrags}
      \psfrag{p1}{\small{$p_1$}}
      \psfrag{p2}{\small{$p_2$}}
      \psfrag{p3}{\small{$p_3$}}
      \psfrag{WUP1}{\small{$W^u(p_1)$}}
      \psfrag{WSP1}{\small{$W^s(p_1)$}}
      \psfrag{WUP2}{\small{$W^u(p_2)$}}
      \psfrag{WSP2}{\small{$W^s(p_2)$}}
      \psfrag{WUP3}{\small{$W^u(p_3)$}}
      \psfrag{WSP3}{\small{$W^s(p_3)$}}
    \includegraphics[width= 2.8in,angle=0]{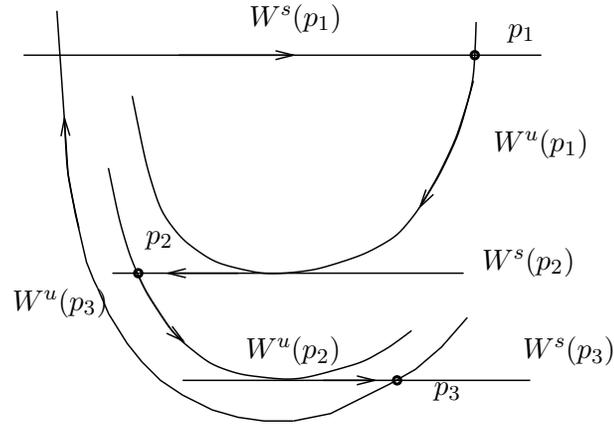}
    \end{psfrags} 
   \caption{ A $3$-floor tower.}
  \end{center}
\end{figure}

{\it{Proof:}}\ \  We prove this Lemma using localized perturbation technic.
As usual consider normal coordinates for a nonresonant saddle $p$.
Induce coordinates in $U$ by normal coordinates for the point $p$ and 
the diffeomorphism $f$.
Application  of Proposition \ref{return} gives existence of
the contour described on Fig.4 in the case $k=3$. Indeed, consider 
an increasing  sequence of numbers $n_1,\dots ,n_k$ such that 
for each $i=1,\dots, k$ the following two properties hold:

1) $T_{n_i}$ intersects $f^{n_i+N}(T_{n_i})$ and they  
form a horseshoe;

2) $n_{i+1}$ is the largest number such that
$T_{n_{i+1}}$ and $f^{n_i+N}(T_{n_i})$ intersect in
a hor\-se\-s\-hoe-like way, i.e., that they bound an open set.

\begin{figure}[htbp]
  \begin{center}
   \begin{psfrags}
      \psfrag{p1}{\small{$p_1$}}
      \psfrag{p2}{\small{$p_2$}}
      \psfrag{p3}{\small{$p_3$}}
      \psfrag{TN1}{\small{$T_{n_1}$}}
      \psfrag{TN2}{\small{$T_{n_2}$}}
      \psfrag{TN3}{\small{$T_{n_3}$}}
    \includegraphics[width= 3in,angle=0]{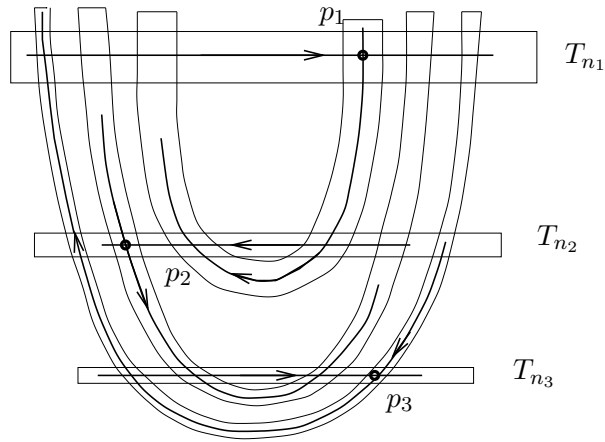}
    \end{psfrags}
   \caption{ An uncomplete $3$-floor tower.}
  \end{center}
\end{figure}

For each $i=1,\dots, k$  condition 1) implies existence of a saddle
periodic point $p_{i} \in T_{n_i} \cap f^{n_i+N}(T_{n_i})$ of period 
$n_i+N$ and condition 2) that $W^s_{loc}(p_{i+1})$ and 
$W^u_{loc}(p_{i})$ intersect.

Let $U$ be equipped with normal coordinates.
Define the maximal distance
in the vertical direction between $W^s_{loc}(p_{i})$ and 
$W^u_{loc}(p_{i})$ as the maximum of distance between  any two points 
$x \in W^s_{loc}(p_{i})$ and $y \in W^u_{loc}(p_{i})$ such that 
$x$ and $y$ have the same $\bar x$-coordinate.
Denote the vertical distance between centers of $T_{n_i}$ and
$T_{n_{i+1}}$ by $t_i$ (see Fig. 5). By calculation in section 
\ref{s2.2} we get $t_i=\mu^{-n_i}-\mu^{-n_{i+1}}$.

\bprop \label{pert} 
If the saddle $p$ having a homoclinic tangency is at least 
$r$-shrinking, then the ratio
$\frac{s_i-t_i}{t_i^r}$ is arbitrarily small for 
each $i=1,\dots ,k-1$.
\enprop

{\it Proof:}\ \  Let us use notations and quantitative estimate obtained in 
Proposition {\ref{return}}. Let $p$ be $\rho$-shrinking, 
$\rho>r$. Recall that the rectangle $T_n$ is 
centered at $(1,\dt_n=\mu^{-n})$ and  has length $ 3c\mu^{-n/2}$ and 
width $\mu^{-3n/2}$. Notice that the width is much less than $\mu^{-n}$, 
the height of center $\mu^{-n}$. Since $p$ is $\rho$-shrinking and 
$n_i$ and $n_{i+1}$ satisfy the conditions  
$\mu^{-n_{i+1}}> const \ \mu^{\rho n_i} >\mu^{-n_{i+1}-1}$
it implies that $s_i-t_i<\dt_{n_{i+1}}+\Delta_{n_i}<C \dt_{n_i}^\rho<
\eps \dt_{n_i}^r=\eps t_{n_i}$ for any $\eps>0$ and a sufficiently
large $n_i$ (see Fig.2 right and Fig.5). Q.E.D. 

\bprop \label{pert1} If the ratio $(s_i-t_i)/t_i^r$ is
arbitrarily small, then there exists a small $C^r$-perturbation 
inside of the ball $B$ (see Fig.5)  such that $W^s_{loc}(p_{n_{i+1}})$
and  $W^u_{loc}(p_{n_i})$ have a point of a heteroclinic tangency.
\enprop
{\it Proof:}\ \ Using the standard perturbation technic
we lift up $W^u_{loc}(p_{n_i})$ and create a heteroclinic tangency. Q.E.D.

\brm \label{gap?}
In order  to construct a $k$-floor tower one needs to
create a heteroclinic tangency
of  $W^s_{loc}(p_{n_{i+1}})$ and  $W^u_{loc}(p_{n_i})$
by a $C^r$-perturbation. We construct it by 
``bending''  $W^u_{loc}(p_{n_i})$. 

Another way to construct it is by fixing the eigenvalue $\mu>1$ and varying 
the other eigenvalue $\lambda<1$ of the saddle $p$ exhibiting homoclinic 
tangency. See Proposition \ref{return}: the rectangle $T_n$ is centered
at $(1,\mu^{-1})$ and the curvilinear rectangle
$f^{n+N}(T_n)$ is $C \dt_n^\rho=C \lambda^n$ distant away from
$W^s(p)$, therefore, by changing $\lambda$ one can vary
the position of $f^{n+N}(T_n)$ without changing the position
of $T_n$. But, in this case one needs some additional geometric 
argument to construct all heteroclinic tangencies of a $k$-tower
simultaneously.
\erm
 
\begin{figure}[htbp]
  \begin{center}
   \begin{psfrags}
    \psfrag{B}{\small{$B$}}     
    \psfrag{ti}{\small{$s_i$}}
    \psfrag{si}{\small{$t_i$}}
    \psfrag{pi}{\small{$p_{2i+1}$}}
    \psfrag{pi+1}{\small{$p_{2i+2}$}}
    \psfrag{dni}{\small{$\dt_{n_i}$}}
    \psfrag{dni+1}{\small{$\dt_{n_{i+1}}$}}
    \psfrag{q}{\small{$q$}}
    \psfrag{WSP}{\small{$W^s(p)$}}
    \psfrag{WUPi}{\small{$W^u(p_{2i+1})$}}
    \psfrag{WSPi}{\small{$W^s(p_{2i+1})$}}
    \psfrag{WUPi+1}{\small{$W^u(p_{2i+2})$}}
    \psfrag{WSPi+1}{\small{$W^s(p_{2i+2})$}}
    \psfrag{app dnir}{\small{$\sim \dt_{n_i}^\rho$}}
    \includegraphics[width= 3in,angle=0]{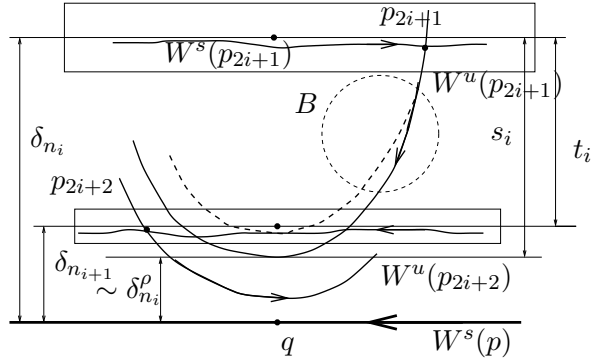}
    \end{psfrags} 
   \caption{A localized perturbation for a floor of a tower.}
  \end{center}
\end{figure}

\bprop \label{shrink} Given $r \in \Bbb Z_+$
if $n_1=n_1(r)$ is sufficiently large, then all saddles
$p_1, \dots, p_k$, described above, are at least $r$-shrinking.  
\enprop 

{\it{Proof:}}\ \ Recall that $p_1, \dots, p_k$ are saddle periodic 
points located in a neighborhood of a homoclinic tangency of a
dissipative saddle $p$. Let $\lb<1<\mu$ denote eigenvalues of $p$. 

With the notations of Lemma \ref{normf}  the return map 
$f^{n+N}$ defined on the rectangle $T_n$ (see section 2.2)  
has the form 

\beq 
\left( \begin{array}{c}
\bar x \\
\bar y
\end{array} \right)
\mapsto
\left( \begin{array}{c}
\alpha \bar y_n+\beta \lb^n (1+\bar x)+ H_1(\cdot,\cdot) \\
\gamma y_n^{2} +\sigma \lb^n (1+\bar x)+ H_2 (\cdot, \cdot)
\end{array} \right),
\eneq
where $\bar y_n=\mu^n \bar y -1,\ 
H_j(\cdot, \cdot)=H_j(\lb^n(1+\bar x),\bar y_n),\ j=1,2$.
The functions $H_1$ and $H_2$ satisfy condition (\ref{remainder}). 
Compare with (\ref{compos}) for $k=2,\ \mu_i=0, i=0,1$.

Therefore, in our notations the fixed point equation has the form

\beq
\begin{cases}
\alpha \bar y_n+\beta \lb^n (1+\bar x)+H_1(\lb^n(1+\bar x), \bar y^n)
=\bar x, \\
\gamma \bar y_n^2+\sigma \lb^n (1+\bar x)+H_2(\lb^n(1+\bar x), \bar y_n)
=\bar y
\end{cases}
\eneq 
Since, the saddle $p$ is dissipative all terms
$\beta \lb^n (1+\bar x),\ H_1(\lb^n(1+\bar x), \bar y^n),\  
\sigma \lb^n (1+\bar x),$ and  $H_2(\lb^n(1+\bar x), \bar y_n)$
tend to zero as $n$ tends to $+\infty$.  
Denote one of fixed points by $(\bar x^0_n,\bar y^0_n)$.
It should belong to $T_n$ which implies that $\bar y^0_n =\mu^{-n} (1+o(1))$.  
Direct calculation of the linear part at $(\bar x^0_n,\bar y^0_n)$
shows that eigenvalues of the linearization approach $2\mu^n$ and $\lb^n/2$ 
respectively. Therefore, if $n$ is sufficiently large, then
$p_n$ is obviously $r$-shrinking for any ahead given $r$. This completes the proof. Q.E.D.

\subsection{The third step: Construction of a $k-th$ order tangency}

We shall prove that by a perturbation of a
$(k+1)$-floor tower one can create a $k$-th order homoclinic tangency. 
Let us start with a 2-nd order tangency and use induction in $k$ then.

\bprop {\cite {GST}} A perturbation of a 3-floor tower 
can create a 2-nd order homoclinic tangency.
\enprop
{\it Proof:} \ \ {\it{Step 1}}. Let us create a 1-st order tangency of 
$W^u_{loc}(p_2)$
and $W^s_{loc}(p_1)$. Start with a 3-tower on Fig.3.
``Push'' $W^u_{loc}(p_2)$ down $W^s_{loc}(p_3)$. Denote by $\gm$ 
the tongue (the part) of $W^u_{loc}(p_2)$ underneath $W^s_{loc}(p_3)$. 
The tongue $\gm$ is in the sector of the saddle hyperbolic point $p_3$,
therefore, under iteration of $f$ $\gm$ will be stretched along $W^u_{loc}(p_3)$ 
and for some $s$ $f^s(\gm) \cap W^s(p_1) \neq \emptyset$. Varying the size 
of the tongue $\gm$ we can create a heteroclinic tangency (see Fig.6.a) with $i=0$).
Denote a point of tangency by $q^*$. Only two parts of $W^u(p_2)$
are depicted on figure 6 a): first part ---  starting part of $W^u(p_2)$ at $p_2$
and second --- image of $\gm$ after a number of iterations under $f$ 
(in above notations $f^s(\gm)$).

Assume that saddle $p_1$ is nonresonant. Then there is normal
coordinates around $p_1$ linearizing $f$. Induce by $f$
normal coordinates in a neighborhood of $U^*$ of $q^*$. In what
follows we shall use these coordinate systems in $U^*$.

\begin{figure}[htbp]
  \begin{center}
   \begin{psfrags}
    \psfrag{q*}{\small{$q^*$}}     
    \psfrag{q1}{\small{$q_1$}}
    \psfrag{q2}{\small{$q_2$}}
    \psfrag{laaa}{\small{$a)$}}
    \psfrag{labb}{\small{$b)$}}
    \psfrag{lacc}{\small{$c)$}}
    \psfrag{ladd}{\small{$d)$}}
    \psfrag{q1=q2}{\small{$q_1=q_2$}}
    \psfrag{gm}{\small{$\gamma$}}
    \psfrag{pi}{\small{$p_{i+1}$}}
    \psfrag{pi+1}{\small{$p_{i+2}$}}
    \psfrag{WSPi}{\footnotesize{$W^s(p_{i+1})$}}
    \psfrag{WUPi}{\footnotesize{$W^u(p_{i+1})$}}
    \psfrag{WSPi+1}{\footnotesize{$W^s(p_{i+2})$}}
    \psfrag{WUPi+1}{\footnotesize{$W^u(p_{i+2})$}}
    \includegraphics[width= 4in,angle=0]{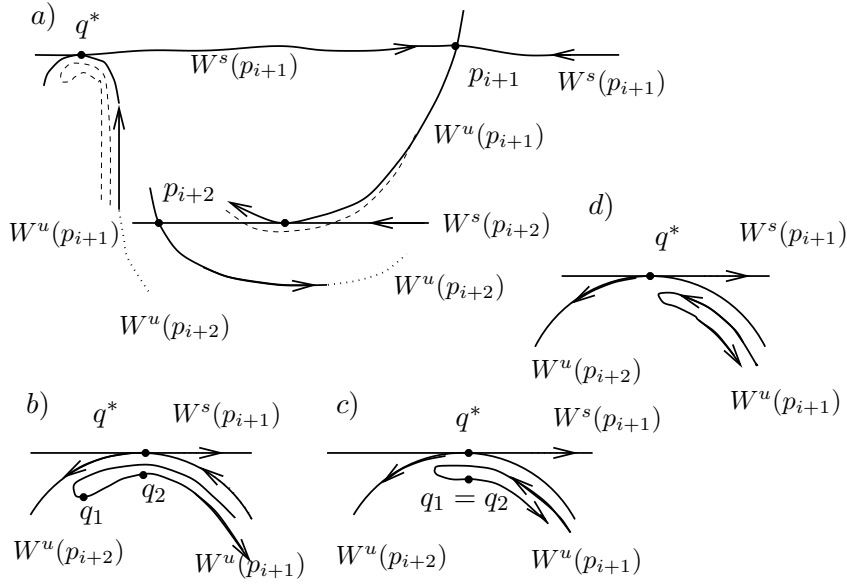}
    \end{psfrags}
   \caption{a 2-nd (even) order tangency.}
  \end{center}
\end{figure}

\begin{figure}[htbp]
  \begin{center}
   \begin{psfrags}
    \psfrag{q*}{\small{$q^*$}}     
    \psfrag{q'}{\small{$\tilde q$}}
    \psfrag{q3}{\small{$q'$}}
    \psfrag{laaa}{\small{$a)$}}
    \psfrag{labb}{\small{$b)$}}
    \psfrag{lacc}{\small{$c)$}}
    \psfrag{gm}{\small{$\gamma$}}
    \psfrag{pi}{\small{$p_{i}$}}
    \psfrag{pi+1}{\small{$p_{i+1}$}}
    \psfrag{WSPi}{\footnotesize{$W^s(p_{i})$}}
    \psfrag{WUPi}{\footnotesize{$W^u(p_{i})$}}
    \psfrag{WSPi+1}{\footnotesize{$W^s(p_{i+1})$}}
    \psfrag{WUPi+1}{\footnotesize{$W^u(p_{i+1})$}}
    \includegraphics[width= 3.5in,angle=0]{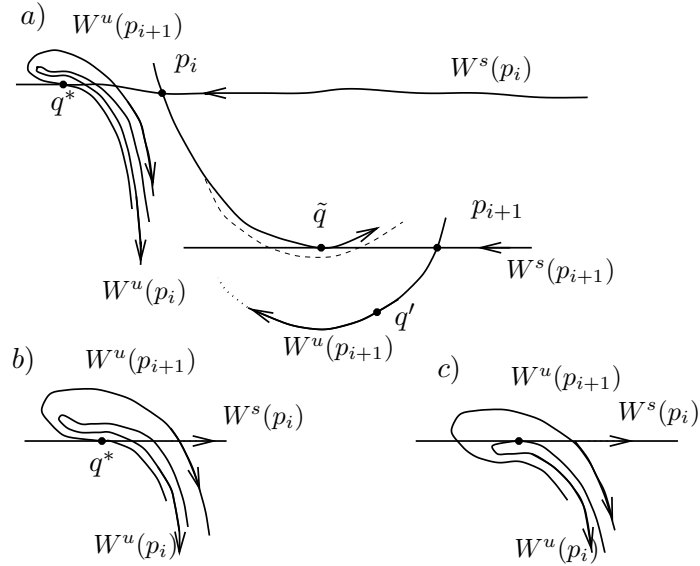}
    \end{psfrags}
   \caption{an odd order tangency.}
  \end{center}
\end{figure}

{\it{Step 2}}. Let us create a 2-nd order homoclinic tangency of 
$W^u_{loc}(p_1)$ and $W^s_{loc}(p_1)$. Start with a contour 
on Fig.6 a). ``Push'' $W^u_{loc}(p_1)$ down $W^s_{loc}(p_2)$. 
Denote by $\gm^1$ the tongue (the part) of $W^u(p_1)$ underneath $W^s(p_2)$.
Some iterate of the tongue $\gm^1$  of $W^u_{loc}(p_1)$ come into $U^*$. 
$U^*$ has normal coordinates and naturally defined the horizontal and the 
vertical directions. Now our goal is varying the size
of $\gm^1$ construct a tangency of some iterate of $\gm^1$ to
the horizontal direction in $U^*$. 

Fix some coordinates in a neighborhood of $\tilde U$ of
a tangency $\tilde q$ of $W^u_{loc}(p_1)$ and $W^s_{loc}(p_2)$.
Consider a $1$-parameter family of diffeomorphisms 
$\{f_\eps\}$, where $\eps$ is the maximal distance of 
$W^u_{loc}(p_1) \cap \tilde U$ and $W^s_{loc}(p_2)\cap \tilde U$ 
in the vertical direction. 

Let $\gm_\eps=W^u_{loc}(p_1) \cap \{y \leq 0\}$. Fix $\eps>0$
and $s=s(\eps$) such that $\gm_\eps^s=f^s(\gm^1_\eps) \cap U^* \neq \emptyset$
(see Fig. 6 b)). The curve $\gm^s_\eps$ has two points
$q_1$ and $q_2$ of tangency to the horizontal direction. As $\eps$
decreases $q_1$ and $q_2$ approach one to the other and for
some critical value $\eps^*$ they collide and $q_1=q_2$
(see Fig. 6 c)).

At the point $q_1=q_2$ $W^u_{loc}(p_1)$ has a 2-nd order
tangency to the horizontal direction. Let this point have
coordinates $(\eps_1,\eps_2)$ in $U^*$. Lifting $W^u_{loc}(p_2)$
by $\eps_2$ we can create a $2$-nd order tangency at point 
$(\eps_1,0)$. This completes the proof of the Proposition.

To construct a $k$-th order tangency assume that we have already 
constructed a diffeomorphism $f \in \textup{Diff}^r(M)$ which has a saddle 
periodic point $p$ exhibiting a homoclinic tangency. In a small
neighborhood of homoclinic tangency there exist two saddle periodic orbits 
$p_1$ and $p_2$ such that $W^s_{loc}(p_1)$ and $W^u(p_2)$ 
has a $(k-1)$-st order tangency and $W^u_{loc}(p_1)$ and $W^s_{loc}(p_2)$
have a tangency. Consider two case $k=2p$ and $k=2p+1$ (see Figures 6 a) 
and 7 a) for $i=0$ respectively). These cases have different topological
pictures. 

\blm \label{lthree} Let $p_1$ and $p_2$ 
be two saddle periodic points and $W^u(p_2)$ have a $(k-1)$-st 
order tangency with $W^s_{loc}(p_1)$ at point $q^*$ and $W^u_{loc}(p_1)$ 
and $W^s_{loc}(p_2)$ have a tangency at a point $\tilde q$. 
Then by a $C^r$-perturbation one can create a $k$-th order 
homoclinic tangency of $W^s(p_1)$ and $W^u(p_1)$
(when $k$ is even see Fig.6 a) for $i=0$ and 
when $k$ is odd  see Fig.7 a) for $i=1$).
\elm
{\it Proof:}\ \  
Assume that $p_1$ and $p_2$ are nonresonant saddles. Fix 
normal coordinate systems $(x_1,y_1)$ (resp. $(x_2,y_2)$) around 
$p_1$ (resp. $p_2$) so that $f$ is linear there.  Let $\lb<1<\mu$ 
be eigenvalues of $p_2$.

Denote by $q'$ a preimage of $q^*$ $q'=f^{-N}(q^*)$.
Fix the normal coordinate systems  $(\hat x,\hat y)$, 
$(\bar x, \bar y)$, and $(\tilde x, \tilde y)$ in neighborhoods
$U^*$ of $q^*$, $\tilde U$ of $\tilde q$, and $U'$ of $q'$
respectively. In what follows we shall use the normal coordinate systems
in $U^*$ and $\tilde U,\ U'$ induced by $f$ from coordinate systems 
$(x_1,y_1)$ and $(x_2,y_2)$ respectively.

The tongue $W^u_{loc}(p_1) \cap \tilde U$ has the form 
$\hat y=a\bar x^2+g(\bar x)$, where $a>0$ and $g(\bar x)=o(\bar x^2)$ at $\bar x=0$.
The map $f^N:\tilde U \to U^*$ has the form
\beq
\left( \begin{array}{c}
\tilde x \\
\tilde y
\end{array} \right)
\mapsto
\left( \begin{array}{c}
\alpha \tilde y_n+ \beta \tilde x + H_1(\tilde x,\tilde y) \\
\gamma \tilde y^{k} +\sigma \tilde x + H_2 (\tilde x, \tilde y)
\end{array} \right),
\eneq
where $H_1(\tilde x,\tilde y)$ and $H_2(\tilde x,\tilde y)$
satisfy condition (\ref{remainder}). 

The idea of the proof is to shift $W^u_{loc}(p_1)$ down to 
$\bar y=a \bar x^2-\eps +g(\bar x)$ and consider the 
versal $k$-parameter family $\{f_\mu\},\ 
\mu=(\mu_0, \dots, \mu_{k-1})$ of the form (\ref{unfold})
unfolding $(k-1)$-st order tangency at the heteroclinic tangency at $q^*$. 
The part of $W^u_{loc}(p_1)$ given by 
$\{\bar y \leq 0\} \cap \{\bar y=a \bar x^2-\eps+g(\bar x)\}$
after a number of iterations under $f$ come to 
a neighborhood $U^*$ of $q^*$. It turns out that by
varying $(k+1)$ parameters $\eps,\mu_0, \dots, \mu_{k-1}$
we can construct a $k$-th order homoclinic tangency in $U^*$.
Let us prove this statement.

Calculate the composition map $f^N_\mu \circ f^n:\tilde U_n \to U^*$,
which is defined in an open subset $\tilde U_n \subset \tilde U$,
\beq \label{comp}
\left( \begin{array}{c}
\hat x \\
\hat y
\end{array} \right)
\stackrel{f^N_\mu \circ f^n}{\longmapsto} 
\left( \begin{array}{c}
\alpha \bar y_n+\ \beta \lambda^n(1+\bar x) +\ 
H_1(\cdot,\cdot)\\
\gamma \bar y_n^{k} +\sum_{i=0}^{k-1} \mu_i \bar y_n^i+
\sigma \lambda^n(1+\bar x)+H_2(\cdot,\cdot) 
\end{array} \right).
\eneq
where $\bar y_n=\mu^n (a\bar x^2 -\eps + g(\bar x)) +1$, 
$H_j(\cdot,\cdot)=\ H_j(\lambda^n(1+\bar x),\bar y_n)$
and $j=1,2$. Assume that after a change of coordinates
in $U^*$ preserving lines $\{\hat y=const \}$ variable $\sigma$ has an 
appropriate sign depending on whether $k$ is odd or even and 
on high order derivatives of $y(x)=a\bar x^2+ g(\bar x)$.

Consider the following parameterization of 
$W^u_{loc}(p_1) \cap \{\bar y \leq 0\}$ by 
$\bar x(t)=t$,\ $\bar y(t)=at^2-\eps +g(t)$.
Let $\hat x_n(t)$ denote the first coordinate function of (\ref{comp})
and $\hat y_n(t)$ --- the second. It is enough to find parameter values
$\eps(n),\ \mu(n)=(\mu_0(n), \dots , \mu_{k-1}(n))$, and
$t^*=t(n)$ such that 

\beq \label{suffcond}
\bar y_n(t^*)=\hat y_n(t^*)=
\left. \frac{\partial \hat y_n(t^*)}{\partial t} \right|_{t=t^*}= \dots=
\left. \frac{\partial^k \hat y_n(t^*)}{\partial t^k} \right|_{t=t^*}=0,
\eneq
provided that $\left.\frac{\partial \hat x_n(t)}{\partial t}\right|_{t=t^*}
\neq 0$. Indeed, $t^*$ corresponds to the point of a
$k$-th order homoclinic tangency of $W^u(p_1)$ and $W^s(p_1)$,
because $f^N_\mu \circ f^n(W^u_{loc}(p_1)) \subset W^u(p_1)$.

In what follows we shall calculate  approximately values of all $k+2$ parameters,
including $t^*$, that  satisfy (\ref{suffcond}). Now we start calculating parameter 
values $t^*$, $\mu_0(n),\ \mu_1(n)$, and so on.

{\it{Step}} 1.
\beq \label{first}
\bar y_n(t)=\mu^n(at^2+g(t)-\eps)+1.
\eneq
Recall that $g(t)$ is $o(t^2)$ at $t=0$. Thus, 
$\bar y_n(t)=0$ for $t^*\approx \sqrt{a^{-1} (\eps-\mu^n)}$.

To simplify calculations notice that in a view of condition 
(\ref{remainder}) $H_2(\lb^n(1+t),0) \sim \lb^{2n}$ for large $n$'s.
By Proposition \ref{shrink} one can choose the saddle $p_1$ to be
at least $(k-1)$-shrinking, so eigenvalues of $p_1$
satisfy $\lb \mu^{k-1}=\tau<1$.  
We shall see that in this case we can choose $\eps$ in such  a way that 
$\mu^n\ t^*(n) \to 0$ as $n \to \infty$
exponentially fast. Without loss of genericity we can choose the 
saddle $p_1$ so that $\lb \ \mu^{k-1}<\tilde \tau<1$ (see Proposition \ref{shrink}). 

Introduce some notations.
$T_n=2a\mu^n \ t^*$,\ 
$C_{sj}(t)=\partial^j \left((\bar y_n(t))^s\right)/ \partial t^j$
for $s,j=1,\dots,k$,\  $\partial^j_r H_2(x^0_1,x^0_2)=
\partial^j H_2(x_1,x_2)/\partial x_r^j|_{x_1=x_1^0,x_2=x_2^0}$.

In order to proceed we need the following 
\bprop \label{der}
There is a set of constants $\{C_{sj}\}_{j \leq s}$,\ $s,j=1,\dots,k$ 
independent of $n$ such that 
\beq \nonumber 
C_{sj}(t^*(n))=
\begin{cases}
(1+o(1))C_{sj}T_n^j \left(t^*\right)^{j-s} \ \ {\textup{for}}\ 
j \leq s \leq 2j \\ 
(1+o(1))C_{sj} T_n^j \left(t^*\right)^{-j} \ \ {\textup{for}}\ s >2j.
\end{cases}
\eneq
\enprop

{\it{Proof of the Proposition:}}\ \ Recall that $\bar y_n(t^*)=0$. So,
it is easy to see that $C_{jj}(t^*)=\left[ \mu^n (2at^* +g'(t^*))\right]^j$ 
$=(1+o(1))T_n^j$.
For $s>j$ one can check that for some positive explicitly computable constant 
$C_{sj}$ we have  $C_{sj}(t^*)=\partial^{s-j}C_{jj}(t)|_{t=t^*}=
(1+o(1))C_{sj}T_n^j\left(t^*\right)^{j-s}$ for $j < s\leq 2j$.
Plugging in the value for $T_n$ we obtain $C_{(2j)j}(t^*)=(1+o(1))C_{(2j)j}\mu^{jn}$. 
For $s>2j$ $C_{sj}(t^*)=(1+o(1))C_{sj}T_n^j \left(t^*\right)^{-j}$.
This completes the proof. Q.E.D.

{\it{Step}} 2.
\beq \label{second}
\hat y_n(t^*)=\mu_0(n)+\sigma \lb^n (1+t^*)+o(\lb^n)=0.
\eneq
Thus, $\mu_0(n) \approx -\sigma \lb^n (1+t^*)$.

Notice that condition (\ref{remainder}) implies 
$\partial^j_2 H_2(\lb^n(1+t^*),\bar y_n(t^*)) \sim \lb^{n}$ for each $j=1,\dots,k$. 
For some $C>0$

\beq
\begin{aligned}
\left| \left. \frac{\partial^j}{\partial t^j} 
H_2(\lb^n(1+t),\bar y_n(t)) \right|_{t=t^*} \right| \leq 
\left| \left. \sum_{s=0}^j \  \partial^{j-s}_2 \ 
\partial^s_1 H_2(\lb^n (1+t),\bar y_n(t)) \right|_{t=t^*}\right| \leq \\
\leq \left| \left. \sum_{s=0}^j \ \lb^{ns}\ \partial^{j-s}_2 \ 
(\partial^s_1 H_2)(\lb^n (1+t),\bar y_n(t)) \right|_{t=t^*}\right|
\leq C \lb^n \mu^{[j/2]n},
\end{aligned}
\eneq 
provided that $T_n \to 0$ as $n\to \infty$. The last
inequality can be proven using formulas from Proposition \ref{der}.

Direct calculation shows that 
because of condition (\ref{remainder}),
$\bar y_n(t^*)=0$, and $T_n \to 0$ as $n\to \infty$ derivative 
$\partial_t H_2(\lb^n(1+t^*),\bar y_n(t^*))=o(\lb^n)$.
 
{\it{Step}} 3.
\beq \label{third}
\left.\frac{\partial \hat y_n(t)}{\partial t}\right|_{t=t^*}=
C_{11}(t^*)\mu_1 +\sigma \lb^n +o(\lb^n)=0.
\eneq
Therefore, $\mu_1(n) \approx - \sigma \lb^n C^{-1}_{11} [T_n]^{-1}$.

{\it{Step}} 4.
\beq \label{forth}
\ \ \quad \left.\frac{\partial^2 \hat y_n(t)}{\partial t^2} \right|_{t=t^*}=
C_{22}(t^*)\mu_2(n)+(1+o(1))C_{21}T_n \left( t^* \right)^{-1}\mu_1(n)+C\lb^n \mu^n=0.
\eneq
By Proposition \ref{der},  $C_{22}(t^*)=(1+o(1))C_{22}T_n^2$. 
Now denote $D_2=C_{21}/C_{22}$.
Thus, $\mu_2(n) \approx -D_{2}\left[T_n t^* \right]^{-1}\mu_1(n)$.

{\it{Step}} 5.
\beq 
\begin{aligned} \label{fifth}
\left.\frac{\partial^3 \hat y_n(t)}{\partial t^3} \right|_{t=t^*}=
C_{33}(t^*)\mu_3(n)+(1+o(1))C_{32}T_n^2 \left( t^* \right)^{-1}\mu_2(n)+\\
+(1+o(1))C_{31}T_n \left( t^* \right)^{-2}\mu_1(n)+C\lb^n \mu^n=0.
\end{aligned}
\eneq
By Proposition \ref{der}, $C_{33}(t^*)=(1+o(1))C_{33}T_n^3$.
Denote $D_3=(C_{32}+C_{31}C_{21})/C_{33}$. Then using the formula
for $C_{33}(t^*)$ we obtain
$\mu_3(n) \approx -1/C_{33}(C_{32}\left[T_n t^* \right]^{-1}\mu_2(n)
-C_{31}\left[T_n t^* \right]^{-2}\mu_1(n))$. Thus, 
$\mu_3(n) \approx- D_3\left[T_n t^* \right]^{-2}\mu_1(n)$.

{\it{Step}} $k+1$.
\beq 
\begin{aligned} \label{(k-1)-st}
\left. \frac{\partial^{k-1} \hat y_n(t)}{\partial t^{k-1}}\right|_{t=t^*}=
C_{(k-1)(k-1)}(t^*) \mu_{k-1} + \mu_0(n)+C \lb^n \mu^{[(k-1)/2]n}\\
\sum_{j=1}^{k-2}(1+o(1))C_{(k-1)(k-j)}T_n^{k-j} 
\left( t^* \right)^{-j}\mu_{k-j}(n)=0.
\end{aligned}
\eneq
By Proposition \ref{der}, we have
$C_{(k-1)(k-1)}(t^*)=(1+o(1))C_{(k-1)(k-1)}T_n^{k-1}$. 
Thus, $\mu_{k-1}(n) \approx -1/C_{(k-1)(k-1)}\sum_{j=0}^{k-2}(1+o(1))C_{(k-1)(k-j-1)} 
\left[T_n t^* \right]^{-j}\mu_{k-j-1}(n)$. For some explicitly 
computable constant $D_{k-1}$,  we obtain
$\mu_{k-1}(n) \approx D_{k-1} \left[T_n t^* \right]^{-k+2}\mu_1(n)$.

At last 
\beq 
\begin{aligned} \label{k-th}
\left. \frac{\partial^{k} \hat y_n(t)}{\partial t^{k}}\right|_{t=t^*}=
C_{kk}(t^*) \gm+
\sum_{j=1}^{k-1}(1+o(1))C_{k(k-j)}T_n^{k-j} 
\left( t^* \right)^{-j}\mu_{k-j}(n)+\\
+\mu_0(n) +C \lb^n \mu^{[k/2]n}=0.
\end{aligned}
\eneq
By Proposition \ref{der} $C_{kk}(t^*)=(1+o(1))C_{kk}T_n^{k}$.
Combining approximate equalities for 
$\mu_j(n)\approx - D_j\left[T_n t^* \right]^{-j+1}\mu_1(n)$ we
obtain
$\gm \approx -1/C_{kk}\sum_j C_{kj} \left[T_n t^* \right]^{-j} \mu_j(n)=
-D_k \left[T_n t^* \right]^{-k+1} \mu_1$ for some explicitly 
computable $D_k$. By a perturbation of the equation $y(x)=ax^2-g(x)$ 
we can guarantee that $D_k$ is different from $0$.

Let us plug in the equation for $\mu_1(n)$.
\beq \label{gamma}
\ \ \quad \gm \approx -\sigma D_k \lb^n T_n^{-1} \left[T_n t^* \right]^{-k+1}=
-(1+o(1))\sigma D_k (2a)^{-k}\lb^n \mu^{-kn} \left(t^* \right)^{-2k+1}.
\eneq
Now depending on the sing of $D_k$, we can choose the sign for $\sigma$ 
above so that both parts of this approximate equality have the same sign. 
By choosing an appropriate $\eps(n)$ we have the following the assymptotic 
formula  $t^* \approx \sqrt{a^{-1}(\eps-\mu^n)} \approx 
\left(\lb \mu^{-k}\right)^{n/(2k-1)}$. Thus, we can satisfy  
assymptotic formula (\ref{gamma}) for $\gm$. 
Let us check the condition $\mu^n t^* \to 0$ exponentially fast.
Since $p_1$ is at least $k-1$ shrinking 
$\mu^n t^*=(\lb \mu^{k-1}))^n=\tau^n<1$. This
complete the proof of Lemma \ref{lthree}. Q.E.D.

To complete the proof of Theorem \ref{Shil} we apply Corollary
\ref{ktang} which allow a diffeomorphism with a $k$-th order tangency
by a $C^r$-perturbation turn into  a diffeomorphism with an arbitrarily 
high period $k$-degenerate orbit.

\bigskip

\subsection{A Proof of Theorem \ref{mthm}}
\label{growth}

Fix a $C^r$ metric $\rho_r$ in Diff$^r(M)$ defined by 
the standard way (see e.g. {\cite{PM}}). 
Let $f$ be a $C^r$ diffeomorphism which 
belongs to a Newhouse domain $\Cal N$.  Write $f \mapsto_{\varepsilon,r} g$
if $g$ is a $C^r$-perturbation of size at most $\varepsilon$ with respect to
$\rho_r$. Consider an arbitrary sequence of positive 
integer numbers $\{a_n\}_{n=1}^\infty$.

Now for any $\varepsilon$ we construct a $3\varepsilon$ 
perturbation $f_3$ of a diffeomorphism $f$ such that for some 
$n_1$ the diffeomorphism $f_3$ 
{\it has $n_1 a_{n_1}$ hyperbolic periodic orbits 
of period  $n_1$}. Hyperbolicity implies that the same is 
true for all  diffeomorphisms sufficiently close to $f_4$.
 
Step 1. $f \mapsto_{\varepsilon,r} f_1$, where 
$f_1$ belongs to a Newhouse domain and is $C^\infty$ smooth.

Step 2. By theorem {\ref {Shil}}, there
exists a $C^r$-perturbation $f_1 \mapsto_{\varepsilon,r} f_2$ such that
$f_2$ has a $k$-degenerate periodic orbit $q$ of an  arbitrarily 
large period, where $k \geq r$.

Step 3. Let $n_1$ be a period of the $k$-degenerate periodic orbit $q$. 
It is easy to show that one can find 
$f_2 \mapsto_{\varepsilon,r} f_3$ such that in a small neighborhood 
of $q$ $f_3$ has $n_1 a_{n_1}$ {\it{hyperbolic 
periodic points}} of period $n_1$. 

Therefore, we show that an arbitrary $C^r$-close to $f$ there exist 
a neighborhood $U \subset {\textup{Diff}}^r(M)$ with the following property
for all $g \in U$
\beq \label{break}
\frac {\#\{x: g^{n}(x)=x\}}{a_{n}} \geq n.
\eneq

If the diffeomorphism $f_1$ belongs to a Newhouse domain 
$\Cal N \subset {\textup{Diff}}^r(M)$, then we can choose perturbation in 
steps 1--3 so small that $f_3$ belongs to the same Newhouse 
domain $\Cal N$. It is not difficult to see from steps 1--3 that 
for an open dense set in $\Cal N$ the condition (\ref{break}) holds 
at least for one $n$. Iterating steps 1--3 one contructs
a residual set such that for each diffeomirphism $f$
from that residual set the condition (\ref{break}) holds
for an infinitely many $n$'s. This completes the proof of 
Theorem \ref{mthm} in the case $\dim M=2$. Q.E.D.

Note that similar inductive argument leads to the well-known
Newhouse's phenomenon on infinitely many coexisting sinks
{\cite {N}}, {\cite {PT}}, {\cite {R}}, and {\cite {TY}}.

\brm {\it{No generic finite parameter family intersects  
a residual set $\Cal R_a$ described in Theorem \ref{mthm} 
with $a_n$ growing quckly enough, e.g. $a_n>n^n$.}}
Indeed, it follows
from the step 4 of previous section, where we take a highly 
degenerate orbit and perturb it. This operation is of large 
codimension and the larger step 
of our induction the more codimension. So, a member of $\Cal R_a$
has to be in an arbitrary small neighborhood of a set of
diffeomorphisms with an arbitrary high degenerate periodic orbit.
\erm

\brm
As we mentioned in the introduction, Theorem \ref{mthm} implies 
Theorem \ref{main}.
\erm

\bigskip
\goodbreak
\subsection{A Proof of Theorem \ref{mthm} in the case $\dim M>2$.}\label{multi}

We shall use the construction described in Step 1 of the
proof of Theorem \ref{mthm}. 

Consider a compact manifold $M$ of dimension $\dim M>2$ and 
a diffeomorphism $F\in \textup{Diff}^r(M)$. Fix a sequence
of numbers $\{a_n\}_{n \in \Z_+}$. Suppose
$F$ has a $C^r$-stable invariant two-dimensional manifold 
$N \subset M$ and the restriction diffeomorphism $f=F|_N:N \to N$ 
belongs to a Newhouse domain $\Cal N \subset {\textup{Diff}}^r(N)$.
$C^r$-{\it{stabililiby of the invariant manifold}} $N$ means that any 
$C^r$-perturbation $\~ F\in \textup{Diff}^r(M)$ of $F$ also has a 
two-dimensional invariant manifold $\~ N$ which is $C^r$-close to $N$ 
and induces a diffeomorphism $\~ f=\~ F|_N:N \to N$  which is $C^r$-close to 
the restriction $f=F|_N:N \to N$ (see Step 1 in \ref{start} below
for an exact formula of $\~ f$). The Fenichel theorem \cite{Fe} gives 
an explicit condition when $F$ has a  $C^r$-stable invariant manifold.
It is important that this is an open condition in $\textup{Diff}^r(M)$.

It is proved in the last subsection that the set of 
diffeomorphims for which the condition (\ref{break}) is satisfied 
for at least one $n\in \Z_+$ is open and dense in a Newhouse domain
$\Cal N \subset \textup{Diff}^r(N)$.
This implies that in a neighbourhood $U$ of $F$ in the space
of diffeomorphisms $\textup{Diff}^r(M)$
there is an open and dense set $\Cal D_1$ of diffeomorphisms 
such that each one satisfies the condition (\ref{break}) for some 
$n\in \Bbb Z_+$.

Let $\Cal D_{1,n_1}$ be an open subset of $\Cal D_1$ consisting
of diffeomorphisms for which the condition (\ref{break}) 
holds (substituting $n=n_1$). There is an open and dense
subset $D_{2,n_1,n_2}$ inside $\Cal D_{1,n_1}$ of diffeomorphisms satisfying 
the condition (\ref{break}) (substituting $n=n_1$ and $n=n_2$ with $n_2>n_1$).
The union  $D_2=\cup_{n_2>n_1} D_{2,n_1,n_2}$ is the open dense set 
inside $U$. Inductive application of these arguments shows that 
there exists a residual set inside of the open set $U$ the condition 
(\ref{break}) holds for infinitely many of $n=n_i$'s. This completes the proof of
Theorem \ref{mthm}. Q.E.D.

\section{A Proof of Theorem \ref{AMthm}}\label{ArMa}
The original proof of Artin-Mazur \cite{AM} uses the fundamental theory of
real algebraic approximations due to Nash. However, this method 
does not work for nonisolated periodic points (see pg.84 {\cite{AM}}).
The method presented below is based on three basic tools:
the Weierstrass approximation theorem, Elimination theory {\cite{Mu}}, and
the Fenichel theorem {\cite{Fe}}. The first two are well-known.
The third is a deep fact about {\it{stability of invariant manifolds}} 
from dynamical systems theory. This method is advantageous 
in that it is simple 
and takes care of nonisolated (even nonhyperbolic) periodic points.

Let us start discribing the proof. 
Consider a $C^r$ diffeomorphism $f:M \to M$. 
We shall approximate $f$ by an A-M diffeomorphism with only
hyperbolic periodic orbits.
There are two steps:

\subsection{Reduction to a problem for polynomial maps}
\label{start}
Using the Whitney Embedding theorem embed $M$ into $\Bbb R^N$ 
for $N=2\dim M+1$. 
Denote by $T$ a tube neighborhood of $M$. For any fixed $r \in \Z_+$
one can extend $f:M \to M$ to a diffeomorphism  $F: T \to T$ such that 
$F$ restricted to $M$ coincides with $f$ and if $F$
contracts along transversal to $M$ directions sufficiently strongly,
then by the Fenichel theorem {\cite{Fe}} each diffeomorphism
$\tilde F:T \to T$ which is $C^r$-close to $F$ has a $C^r$ smooth
invariant manifold $\tilde M$ which is $C^r$-close to $M$. 
Denote by $\pi: \tilde M \to M$ a diffeomorphism from $\tilde M$ to $M$
which can be obtained by projection along the normal to $M$ directions.
Then $\tilde f=\pi^{-1}\circ\tilde F|_{\tilde M}\circ\pi: M \to M$ is a
diffeomorphism which is $C^r$-close to $f$. 
By the Weierstrass approximation theorem one can approximate
a diffeomorphism $F:T\to T$ of an open set $T$ in the Euclidean space
$\R^N$ into itself by a polynomial diffeomorphism
$\tilde F=P|_T:T \to T$. Notice that if $\tilde F$ has only hyperbolic 
periodic orbits, then the induced diffeomorphism 
$\tilde f=\pi^{-1}\circ\tilde F|_{\tilde M}\circ\pi: M \to M$
also has only hyperbolic periodic orbits.

We shall prove that, indeed, one can approximate any
diffeomorphism $F:T\to T$ by a polynomial diffeomorphism
$\tilde F=P|_T:T \to T$ which has only hyperbolic periodic orbits.

Let $D \in \Z_+$. Denote by $A_N^D$ the space of  
vector-polynomials $P:\Bbb R^N \to \Bbb R^N$ of degree 
at most $D$. If $\mu=\mu(N,D)=\#\{\al \in \Z^N_+: |\al| \leq D\}$, 
then $A_N^D$ is isomorphic to $\R^\mu$. Consider $A^D_N$ with 
the Lebesgue measure on it. In what follows we call a 
{\it{vector-polynomial}} by a 
{\it{polynomial}} for brevity.
 
Step 2. {\it{For any $D \in \Bbb Z_+$ an almost every 
polynomial $P:\Bbb R^N \to \Bbb R^N$ from $A_N^D$
has  only hyperbolic periodic orbits and their number  
grows at most exponentially}}. 

The second part of this statement is easy 
provided that the first is true. Indeed, fix 
$k \in \Bbb Z_+,\ k>0$ and consider the system
$$
P(x_1)-x_2=0,\ P(x_2)-x_3=0, \dots,\ P(x_k)-x_1=0.
$$ 
This system has $Nk$ equations, each of them of degree at most $D$. 
By the Bezout theorem  the number of isolated solutions is at most 
$D^{kN} \leq (D^N)^k.$ If all periodic points are hyperbolic,
then they are all isolated and this completes the proof.

Fix $k \in \Bbb Z_+,\ k>0$. Let $\alpha=(\al_1,\dots,\al_N)
\in \Z^N_+$
be a multiindex, $|\al|=\sum_i \al_i$. Fix a coordinate
system in $\R^N$ so one can write each polynomial 
$P(a,\cdot):\Bbb R^N \to \Bbb R^N$  from $A_N^D$ in the form
\beq \label{notat}
\begin{aligned}
P(a,x)=\sum_{|\al |\leq D} a_\al x^\al,&\ \ \textup{where}\  
\ a=(\{a_\al\}_{|\al |\leq D}) \in \R^{\mu},
\\ 
\ x=(x_1,\dots ,x_N) \in \R^N,&\ \ \textup{and}\ \  
x^\al=x_1^{\al_1} \dots x_N^{\al_N}.
\end{aligned}
\eneq 

\blm \label{artmaz} Let $\lb_0 \in \C$ and $|\lb_0|=1$. 
For any $D \in \Bbb Z_+$ an almost every 
polynomial $P:\Bbb R^N \to \Bbb R^N$ from $A_N^D$
has no periodic orbits with the eigenvalue $\lb_0$.
\elm

Denote $P^{(k)}(a,\cdot)=P(a,\cdot) \circ \dots 
\circ P(a,\cdot):\R^N \to\R^N$ ($k$ times), the linearization 
matrix of the map $P^{(k)}(a,\cdot)$ at a point $x$ by 
$d_x (P^{(k)} )(a,x)$, and the $N\times  N$ identity matrix
by $Id$. Let $\lb \in \C$ be a complex number. Denote 
$D(a,\lb,x)=\det\ (d_x\left(P^{(k)}\right)(a,x)-\lb \ Id)$.
Every periodic orbit of period $k$, which has an eigenvalue $\lb$, 
satisfies the following system:

\beq \label{hyperb}
\begin{cases}
P^{(k)} (a,x)-x=0,\  \ x=(x_1, ... x_N) \in \Bbb R^N\\
D(a,\lb,x)=0,\ \  a \in \R^\mu 
\end{cases}
\eneq 
The general goal is to prove that for a ``generic'' choice of 
coefficients $a \in \R^\mu$ of $P(a,\cdot)$ this system
has no solutions satisfying the condition $|\lb|=1$ or 
there is no nonhyperbolic periodic orbit of period $k$.
First, we prove that  a ``generic'' choice of 
coefficients $a \in \R^\mu$ of $P(a,\cdot) $ has no periodic
points with the eigenvalue $\lb=\lb_0$.

Notice that the system (\ref{hyperb}) including the condition 
$\lb=\lb_0$ (or $|\lb|=1$) consists of $N+2$ equations and for each 
value $a$ only $N+1$ indetermine variables $x_1,\dots ,x_N,\lb$.
It might be clear intuitively that for a ``generic'' $a \in \R^\mu$
there is no solution, because the number of equations is more
than the number of indeterminates. To prove it rigorously for $\lb=\lb_0$
we shall apply Elimination theory. 

\subsection{Elimination theory}

Let $\C^m$ denote the $m$-dimensional complex space 
$z=(z_1,\dots,z_m) \in \C^m,\ \ m \in \Z_+$.
A set $V$ in $\C^m$ is called {\it{a closed algebraic set}} 
in $\C^m$ if there is a finite set of polynomials 
$F_1,\dots F_s$  in $z_1,\dots,z_m$ such that
\beq \nonumber
V(F_1,\dots ,F_s)=\{(z_1, \dots ,z_m) \in \C^m |\ 
F_j(z_1,\dots,z_m)=0,\ 1\leq j \leq s\}.
\eneq
One can define a topology in $\C^m$, called the 
{\it{Zariski topology}}, whose closed sets are closed
algebraic sets in $\C^m$. This, indeed, defines a topology,
because the set of closed algebraic sets is closed under a
finite union and an arbitrary intersection. 
Sometimes, closed algebraic sets are also called
Zariski closed sets. 

\bdef A subset $S$ of $\C^m$ is called constructible
if it is in the Boolean algebra generated by the 
closed algebraic sets; or equivalently if $S$ is a disjoint
union $T_1 \cup \dots \cup T_k$, where $T_i$ is locally
closed, i.e. $T_i=T'_i-T''_i$, $T'_i$ --- a closed algebraic 
set and $T''_i\subset T'_i$ --- a smaller closed algebraic.
\endef

One of the main results of Elimination theory is 
the following 

\bthm {\cite{Mu}}\label{constr} Let 
$V \subset \C^\mu \times \C^N$ be 
a constructible set and $\pi:\C^\mu \times \C^N \to \C^\mu$
be the natural projection. Then  $\pi(V)\subset \C^\mu$ is a 
constructible set. 
\ethm 

\brm
An elementary description of elimination theory can be found 
in books Jacobson {\cite{J}} and van der Waerden {\cite{W}}.
\erm

\subsection{Proof of Lemma \ref{artmaz} or
application of Elimination theory to the system (\ref{hyperb})}
Put $\lb=\lb_0$  and consider the system (\ref{hyperb}) as if it is 
defined for $(a;x) \in \C^{\mu} \times \C^{N}$. Then
it defines 
a closed algebraic set $V_k(\lb_0) \subset \C^{\mu} \times \C^{N}$. 
By Theorem \ref{constr}  the natural projection 
$\pi:\C^{\mu} \times \C^{N} \to \C^{\mu}$
of $V_k(\lb_0)$, namely $\pi(V_k(\lb_0)) \subset \C^{\mu}$, 
is a constructible set. The only thing left to show is that 
$\pi(V_k(\lb_0)) \neq \C^\mu$ and has a positive codimension.
Recall that $|\lb_0|=1$ and $\C^\mu$ is the space of
coefficients of polynomial of degree $D$.

\bprop \label{example} Let $\R^\mu$ be naturally embedded 
into $\C^\mu$. Then there is an open set $U \subset \C^\mu$
such that $U \cap \R^\mu \neq \emptyset$ and
for any $a \in U$ the corresponding polynomial
$P(a,\cdot):\C^N \to \C^N$ of degree has exactly $D^{Nk}$
periodic points of period $k$ and all of them
are hyperbolic.
\enprop
{\it Proof:}\ \  Consider the homogeneous polynomial 
$P(a^*,\cdot): \C^N \to \C^N$ of degree $D$
$P(a^*,\cdot):(z_1,\dots,z_N) \mapsto (z_1^D\dots ,z_N^D).$
It is easy to see that $P$ has exactly $D^{Nk}$
periodic points of period $k$ all of which are hyperbolic.
From one side hyperbolicity of periodic points of
period $k$ of $P$ implies that any polynomial 
mapping $\~ P$, which is a small
perturbation of $P$, has at least $D^{Nk}$ hyperbolic points 
of period $k$, but from the other side
Bezout's Theorem implies that a polynomial of degree
$D$ has at most $D^{Nk}$ periodic point of period $k$.
Thus, there is a neighborhood $U \subset \C^\mu$
of $a^*$ 
such that for any $a \in U$ the corresponding polynomial
$P(a,\cdot):\C^N \to \C^N$  has only hyperbolic
periodic points of period $k$ and by definition 
$U \cap \pi(V_k(\lb_0))=\emptyset$. 
Since, $\pi(V_k(\lb_0))$ is constructible, this implies that
$\pi(V_k(\lb_0))$ has positive codimension in $\C^\mu$.
This completes the proof of Proposition.

By Proposition \ref{example} the restriction 
$\pi(V_k(\lb_0)) \cap \R^\mu$ has
positive codimension and, therefore, measure zero 
in $\R^\mu$. Thus, almost every polynomial $P(a,\cdot)$ 
from $A^D_N=\R^\mu$ has no periodic points of period $k$ 
with the eigenvalue $\lb_0$.
Intersection over all $k \in \Z_+$ gives that the 
same is true for all periods. This completes the 
proof of Lemma \ref{artmaz}. Q.E.D.

\subsection{Completion of Proof of Step 2 of Theorem \ref{AMthm}}

Consider the system (\ref{hyperb}) as if it is 
defined for $(a,\lb;x) \in \C^{\mu} \times \C \times \C^{N}$. 
It defines a closed algebraic set, denoted by 
$V_k \subset \C^{\mu} \times \C \times \C^{N}$. 
By Theorem \ref{constr}  the natural projection 
$\pi:\C^{\mu} \times \C \times \C ^{N} \to \C^\mu \times \C$
of $V_k$, namely, $S_k=\pi(V_k) \subset \C^{\mu}\times \C$, 
is a constructible set.

Consider natural projections
$\pi_1:\C^\mu \times \C \to \C^\mu$ and 
$\pi_2:\C^\mu \times \C \to \C$. It 
follows from Proposition \ref{example} 
that $S_k$ has dimension $\mu$. Indeed, the projection
$\pi_1(S_k)=W_k$ contains an open set $U \subset \C^\mu$ and
$S_k$ does not intersect a neighborhood of 
$U \times \{\lb:|\lb|=1\} \subset \C^\mu \times \C$.

By Theorem \ref{constr} that $\pi_1(S_k)=W_k$ is 
constructible and by has Proposition \ref{example} 
$\dim W_k=\mu$.
By Sard's lemma for algebraic sets {\cite{Mu}}
there exists a proper algebraic set 
$\Sigma_k \subset S_k$ such that $\~ S_k=S_k \setminus \Sigma$
outside $\Sigma_k$  the restricted to $\~ S_k$ map 
$p_k=\pi_1|_{\~ S_k}:\~ S_k \to \C^\mu$ has no critical points. 
Thus, outside of some proper closed algebraic set 
$\Sigma'_k \subset \C^\mu$ 
the map $p_k:\~ S_k \to \C^\mu \setminus \Sigma'_k$
is locally invertible.

Recall that our goal is to show that 
$Z=\pi_1(S_k\ \cap \{\lb:|\lb|=1\}) \cap \R^\mu$
has measure zero in $A^D_N=\R^\mu$. 
It is sufficient to show that this locally.

Let $a \in \R^\mu \setminus \Sigma'_k$ and 
$U \subset \C^\mu \setminus \Sigma'_k$ be a neighborhood of $a$.
By construction the map 
$p_1:\~ S_k \to \C^\mu \setminus \Sigma'_k$ is locally invertible,
so the preimage $p_1^{-1}(U)$ consists of a finite disjoint union of
open sets $\cup_{j\in J}U_j \subset S_k$. Thus, one can define a finite 
collection of analytic functions
$\{\lb_{k,j}=\pi_2 \circ p_{k,j}^{-1}:U\to \C\}_{j \in J}$, where
$p_{k,j}^{-1}:U \to U_j$ is the inverse of the restriction
$p_k|_{U_j}:U_j \to U$. 
We need to show that 
$$
\cup_{j \in J} \{\lb^{-1}_{k,j}(\{\lb:|\lb|=1\}) \cap \R^\mu
$$
has measure zero. If for some $j \in J$ the function
$\lb_{k,j}:U \to \C$ is equal to a constant $\lb$,
then by  Proposition {\ref{example}} we have $|\lb|\neq 1$ and
 the preimage $\lb_j^{-1}(Z)=\emptyset$. 
If for some $j \in J$ the function
$\lb_{k,j}:U \to \C$ is not constant, then the set
$\lb_j^{-1}(\{\lb:\ |\lb|=1\}) \cap \R^\mu$ is real analytic set.
Moreover, it follows from Proposition \ref{example} and 
the identity theorem that $\lb_j^{-1}(\{\lb:\ |\lb|=1\}) \cap \R^\mu$
has to be nowhere dense and, therefore, have a positive
codimension. A real analytic set of positive codimension 
has measure zero. It follows e.g. from the fact 
that a real analytic can be stratified (see e.g. {\cite{H}} 
or {\cite{GM}}), i.e., in particular, can be decomposed into at 
most countable union of semianalytic manifolds. Each 
semianalytic manifold must have a positive codimension and, 
therefore, measure zero. 
This implies that for almost every 
$a \in \R^\mu$ the system (\ref{hyperb}) has no solutions 
for $|\lb|=1$. Intersection of all periods $k\in \Z_+$ 
completes the proof of step 2.

Let us complete the proof of Theorem \ref{AMthm}. 
Application of Step 1 shows that 
a diffeomorphism $f:M\to M$ can be extended to a tube
neighborhood $T$ of $M$ $F:T\to T$ and that it is sufficient to 
approximate $F$ by a diffeomorphism $\~F:T \to T$ which 
has only hyperbolic periodic points. 
By the Weierstrass approximation theorem $F$ can be 
approximated by a polynomial diffeomorphism $\~F=P|_T:T \to T$.
Since, in the space of polynomial of any degree $D$ 
polynomial maps with only hyperbolic periodic points form
a full measure set one can choose $\~F=P|_T:T \to T$
which has only hyperbolic periodic points. If $\~ F:T \to T$
has only hyperbolic periodic points, then its restriction
to an invariant manifold also has only hyperbolic periodic points.
This completes the proof of Theorem \ref{AMthm}.

\brm
In order to give a positive answer to the Artin-Mazur
question stated in the introduction it is sufficient to
use only Step 1 and Lemma {\ref{artmaz}} of the above
proof. 
\erm

\section{Prevalence and open problems}\label{prev}

There are two point of view on a notion of genericity
in dynamical systems, singularity theory and etc.: topological and metrical.  
Topological genericity which goes back 
to Baire is standard and widely accepted. 
It says that a property of dynamical systems is
generic if systems with that property form a residual set. 

However, it is easy to construct a residual set in the
segment $[0,1]$ which has measure zero.
Different examples from the KAM theory, small denominators, 
fractal geometry, and so on show that topological
description is not always a good one (see {\cite {HSY}} for more 
examples).
 
Let us describe another point of view which goes back to
Kolmogorov. In his
plenary talk on the International Mathematical Congress
in 1954, A.N. Kolmogorov 
proposed to judge whether a phenomenon is generic or not by
considering a generic finite parameter family with the Lebesgue
measure on a parameter space and looking at measure of parameters 
corresponding to that phenomenon.
In {\cite{K}}
the author proposed the following definition:

Let $B^n$ be an $n$-dimensional ball.
Denote by $\textup{Diff}^r(M,B^n)$ the space of $n$-parameters
families of diffeomorphisms $\{f_\varepsilon\}_{\varepsilon \in B^n}$
with the uniform $C^k$ topology.

\bdef We call a set $P \subset \textup{Diff}^r(M)$  an 
$n$-prevalent with respect to an $n$-parameter family  
$\{f_\varepsilon\}_{\varepsilon \in B^n}$ if 
$P$ restricted to that family form a set of full-measure
with respect to the natural Lebesgue measure in the space of 
parameters:
\beq \label{preval}
\textup {mes}\{ \varepsilon \in B^n:
f_\varepsilon \in P \} - {\textup {full-measure}}.
\eneq

We call a set $P \subset \textup{Diff}^r(M)$ a strictly $n$-prevalent if
the following two conditions hold:

A) $P$ is prevalent with respect to  an open dense set of
$n$-parameter families $\{f_\varepsilon\}_{\varepsilon \in B^n}
\in \textup{Diff}^r(M,B^n)$;

B) For any e\-le\-ment $f \in {\textup{Diff}}^r(M)$ there exists
an $n$-parameter family of diffeomorphisms
$\{f_\varepsilon\}_{\varepsilon \in B^n}
\in {\textup{Diff}}^r(M,B^n)$ which passes through $f$, i.e. $f=f_0$ and
$P$ is prevalent with respect to $\{f_\varepsilon\}_{\varepsilon \in B^n}$.

A set  $P \subset \textup{Diff}^r(M)$ is called $n$-prevalent if
it contains a countable intersection of $n$-prevalent sets.
We also call a set prevalent if it is $n$-prevalent for some $n$ and
neglectable if the complement is prevalent.
\endef   
It easily follows from the definition that a countable intersection of 
$n$-prevalent sets is $n$-prevalent.
\bprop {\cite{K}}
If $P \subset \Bbb R^N$ and $P$ is $n$-prevalent for some $n<N$, 
then $P$ has full measure in $\Bbb R^N$.
\enprop
It shows that on the contrary to the topological genericity
test of this definition in a finite-dimensional case gives a satisfactory
result.

In {\cite {K}} it is proven that certain fundamental facts
from the singularity theory and the theory of dynamical systems
such as transversality theorems, the Whitney embedding, the 
Mather stability, and the Kupka-Smale theorems, which are 
topologically generic, are also prevalent.  In {\cite{HSY}} a definition 
of prevalent set in an infinite dimensional linear space is proposed. 

In a view of this it is natural to pose following problems:

{\bf{Problem 1.}}\ 
{\it {Do Artin-Mazur diffeomorphisms form a prevalent set?}}

{\bf{Problem 2.}}\  {\it Whether or not Newhouse's phenomenon on 
infinitely many coexisting sinks $1$-prevalent?}
\footnote{A partial answer is in {\cite {TY}}.}

Recall also a growth problem for vector fields from 
Artin-Mazur's paper {\cite{AM}}:

{\bf{Problem 3.}}\ 
{\it{Let $X$ be a differentiable vector field on a compact 
manifold $M$. Denote by $N_t(X)$ the number of periodic orbits of $X$, 
period less than or equal to $t$. Does $N_t(X)$ grow at most exponentially 
for some reasonable dense class of vector fields?}} 

{\it Acknowledgments:}\ \ 
I would like to express my warmest thanks to my thesis 
advisor John Mather. He proposed to look at the 
problem of growth of number of periodic orbits and 
suggested the pertinent idea that a highly degenerate periodic orbit 
can generate a lot of periodic orbits in an open way. Several 
discussions of elimination theory with him were fruitful
for me. J. Milnor gave the reference to the book of Bowen and 
proposed the excellent title. A moral support of J.Milnor 
was very important for me. A.Katok pointed out to me the question 
posed by Artin - Mazur (see Theorem \ref{AMthm}).
Remarks of D. Dolgopiat, Yu. Ilyashenko, J. Mather, 
J. Milnor, S. Patinkin, and G. Yuan on the text are highly appreciated.  
G. Levin mentioned to me the elimination theory. Let me express sincere 
gratitude to all of them.

\end{document}